\documentclass[a4,11pt]{amsart}
\usepackage{amssymb,amsmath,amsthm,txfonts}
\usepackage{cite}

\usepackage{color}

\newtheorem{thm}{Theorem}[section]
\newtheorem{lem}[thm]{Lemma}

\newtheorem{prop}[thm]{Proposition}
\theoremstyle{definition}

\theoremstyle{remark}
\newtheorem{rem}[thm]{\textbf{Remark}}
\newtheorem{rems}[thm]{\textbf{Remarks}}

 \makeatletter
    
    \@addtoreset{equation}{section}
  \makeatother

\makeatletter
      \def\@makefnmark{%
         \leavevmode
            \raise.9ex\hbox{\check@mathfonts
                \fontsize\sf@size\z@\normalfont%
                            \@thefnmark}%
       }
      \makeatother

\newcommand{\D}{\textrm{div}}

\newcommand{\dd}{\textrm{d}}

\begin{document}

\title[]{The vorticity equations in a half plane with measures as initial data}
\author[]{K.Abe}
\date{}
\address[K. ABE]{Department of Mathematics, Graduate School of Science, Osaka City University, 3-3-138 Sugimoto, Sumiyoshi-ku Osaka, 558-8585, Japan}
\email{kabe@sci.osaka-cu.ac.jp}

\date{\today}

\maketitle

\begin{abstract}
We consider the two-dimensional Navier-Stokes equations subject to the Dirichlet boundary condition in a half plane for initial vorticity with finite measures. We study local well-posedness of the associated vorticity equations for measures with a small pure point part and global well-posedness for measures with a small total variation. Our construction is based on an $L^{1}$-estimate of a solution operator for the vorticity equations associated with the Stokes equations. 
 \end{abstract}

%contents
%\tableofcontents

%セクション1
\section{Introduction}

\vspace{10pt}

We consider the Navier-Stokes equations in a half plane:

\begin{equation*}
\begin{aligned}
\partial_t u-\Delta{u}+u\cdot \nabla u+\nabla{p}= 0,\quad \D\ u&=0  \qquad \textrm{in}\ \mathbb{R}^{2}_{+}\times (0,\infty),  \\
u &=0\qquad \textrm{on}\ \partial\mathbb{R}^{2}_{+}\times (0,\infty), \\
u&=u_0\hspace{18pt} \textrm{on}\ \mathbb{R}^{2}_{+}\times\{t=0\},
\end{aligned}
\tag{1.1}
\end{equation*}\\
for initial data $u_0={}^{t}(u_{0}^{1},u_{0}^{2})\in L^{2,\infty}_{\sigma}(\mathbb{R}^{2}_{+})$ with a finite measure $\omega_0=\partial_1u^{2}_0-\partial_2u^{1}_0\in M(\overline{\mathbb{R}^{2}_{+}})$, where

\begin{align*}
L^{2,\infty}_{\sigma}(\mathbb{R}^{2}_{+})
=\left\{ f\in L^{2,\infty}(\mathbb{R}^{2}_{+})\ \middle|\ \D\ f=0\ \textrm{in}\ \mathbb{R}^{2}_{+},\ f^{2}(x_1,0)=0,\ x_1\in \mathbb{R}      \right\},
\end{align*}\\
and $M(\overline{\mathbb{R}^{2}_{+}})$ denotes the space of finite real regular Borel measures on $\overline{\mathbb{R}^{2}_{+}}$ equipped with the total variation $||\cdot ||_{M}$. Examples of such $\omega_0$ are vortex sheets and point sources of vorticities.  A vortex sheet is a continuous measure supported on a smooth curve in the plane and a point source is a pure point measure. For the Cauchy problem, global-in-time solutions exist for such initial data \cite{Cottet}, \cite{GMO} (see also \cite{Ben94}, \cite{Brezis94}, \cite{Kato94}), while for a half plane a few results is in known. As in $\mathbb{R}^{2}$, initial velocity $u_0\in L^{2,\infty}_{\sigma}$ satisfying $\omega_0\in M$ is represented by the Biot-Savart law 

\begin{align*}
u_0(x)&=\int_{\mathbb{R}^{2}_{+}} \nabla^{\perp}_{x}D(x,y) \omega_0(\dd y),  \tag{1.2}
\end{align*}\\
where $\nabla^{\perp}={}^{t}(\partial_2,-\partial_1)$ and 

\begin{align*}
D(x,y)&=E(x-y)-E(x-y^{*}),\quad y^{*}={}^{t}(y_1,-y_2), \\
E(x)&=-\frac{1}{2\pi}\log{|x|}.
\end{align*}\\
The right-hand side of (1.2) is an integral by the Borel measure $\omega_0$. We write (1.2) by $u_0=K\omega_0$. Since $K\omega_0=\nabla^{\perp}E*{\omega}_0^{\textrm{odd}}$ for a measure ${\omega}_0^{\textrm{odd}}$ and the convolution $*$ in $\mathbb{R}^{2}$, $K$ acts as a bounded operator from $M$ to $L^{2,\infty}$. If the total variation of $\omega_0$ is small, $u_0$ is small in $L^{2,\infty}$. Hence for small $\omega_0\in M$ unique global-in-time solutions to (1.1) exist by a small data result in $L^{2,\infty}$ \cite{KY}. If the total variation of $\omega_0$ is large, even local well-poseness of (1.1) is unknown in general.

We study the vorticity equations associated with (1.1):

\begin{equation*}
\begin{aligned}
\partial_t \omega-\Delta\omega+u\cdot \nabla \omega&= 0 \hspace{34pt} \textrm{in}\ \mathbb{R}^{2}_{+}\times (0,\infty),  \\
\partial_2\omega-A \omega&=-\partial_1 p_1\quad \textrm{on}\ \partial\mathbb{R}^{2}_{+}\times (0,\infty),  \\
\omega&=\omega_0\hspace{27pt} \textrm{on}\ \mathbb{R}^{2}_{+}\times \{t=0\},
\end{aligned}
\tag{1.3}
\end{equation*}\\
where $A$ is the generator of the Poisson semigroup 

\begin{align*}
e^{sA}g(x_1)=\int_{\mathbb{R}}P_s(x_1-y_1) g(y_1)\dd y_1,\quad P_s(x_1)=\frac{s}{\pi(|x_1|^{2}+s^{2})}, 
\end{align*}\\
and $p_1=p-p_2$ is a remainder from the harmonic pressure $p_2=-\int_{x_2}^{\infty}e^{sA}\partial_1\omega\dd s$. By the Fourier transform, we write $A=-H\partial_1$ with the Hilbert transform $H$ (see Section 3 for the definition of $H$). Since $-\Delta u=\nabla^{\perp}\omega$ and $H^{2}=-I$, the boundary condition $(1.3)_2$ follows by taking the tangential trace to $(1.1)_1$. 

The vorticity equations (1.3) is studied in \cite{Maekawa13} by using a solution formula for that associated with the Stokes flow (i.e., $u=0, p_1=0$ in (1.3)), 

\begin{equation*}
\begin{aligned}
T_{0}(t)\omega_0&=\int_{\mathbb{R}^{2}_{+}} W_{0}(x,y,t)\omega_0(y)\dd y,\\
W_{0}(x,y,t)&=\Gamma(x-y,t)-\Gamma(x-y^{*},t)+2(H\partial_1-\partial_2)\partial_2E*\Gamma(x-y^{*},t), \\
\Gamma(x,t)&=\frac{1}{4\pi t}e^{-|x|^{2}/4t}.
\end{aligned}
\tag{1.4}
\end{equation*}\\
The formula (1.4) is written with $A$ in \cite{Maekawa13}. We write it with $H$. Since the Hilbert transform is bounded on $L^{q}(\mathbb{R})$ for $q\in (1,\infty)$, $T_0(t)$ is a bounded operator on $L^{q}(\mathbb{R}^{2}_{+})$. For $q=1$, the kernel $W_0(x,y,t)$ is not integrable in $\mathbb{R}^{2}_{+}$ for the $x$-variable (see Remarks A.2 (iii)). The formula (1.4) is available to represent vorticity of the Stokes flow, provided that the tangential trace of  $u_0$ vanishes, i.e., 

\begin{align*}
u^{1}_{0}(x_1,0)=0,\quad x_1\in \mathbb{R}.  \tag{1.5}
\end{align*}\\
By (1.2), this condition is equivalent to

\begin{align*}
\int_{\mathbb{R}^{2}_{+}}\frac{y_2}{|x_1-y_1|^{2}+y_2^{2}}\omega_0(y)\dd y=0,\quad x_1\in \mathbb{R}.   \tag{1.6}
\end{align*}\\
For example, for $\omega_0\in L^{q}$, $q\in (1,2)$, satisfying (1.6), the harmonic function 

\begin{align*}
\int_{\mathbb{R}^{2}_{+}} \partial_{x_2}E(x-y^{*})\omega_0(y)\dd y
\end{align*}\\
vanishes by the Liouville theorem. Hence, 

\begin{align*}
T_0(t)\omega_0
\to \omega_0+2(H\partial_{x_2}-\partial_{x_1})\int_{\mathbb{R}^{2}_{+}}\partial_{x_2} E(x-y^{*})\omega_0(y)\dd y
=\omega_0\quad \textrm{in}\ L^{q}\quad \textrm{as}\ t\to0.
\end{align*}\\
(We give a proof for (1.4) in Appendix A for the completeness.) If $\omega_0$ is integrable, (1.5) implies the zero total mass condition,  

\begin{align*}
\int_{\mathbb{R}^{2}_{+}}\omega_0(y)\dd y=0,  \tag{1.7}
\end{align*}\\
by integrating (1.6) by the $x_1$-variable \cite{Maekawa13}. 

The condition (1.5) is not always satisfied for all $\omega_0\in M$. For example, if $\omega_0$ is a point mass, e.g., $\omega_0=\kappa \delta_{x_0}$ for $\kappa \in \mathbb{R}$ and the Dirac measure $\delta_{x_0}$ at $x_0={}^{t}(0,1)$, the tangential trace of $u_0=\kappa \nabla^{\perp}D(x,x_0)$ does not vanish, i.e., 

\begin{align*}
u^{1}_{0}(x_1,0)=\frac{\kappa}{\pi(|x_1|^{2}+1)},\quad x_1\in \mathbb{R}.
\end{align*}\\
For $\omega_0\in M$, the tangential trace $u^{1}_0(x_1,0)$ belongs to $L^{1}(\mathbb{R})$ by (1.2). To study (1.3) for measures $\omega_0\in M$, we construct a different solution operator based on the Green matrix of the Stokes semigroup \cite{Sl03}. As is well known, the integral form of (1.1) is 

\begin{align*}
u(t)=S(t)u_0-\int_{0}^{t}S(t-s)\mathbb{P}(u\cdot \nabla u)\dd s,  \tag{1.8}
\end{align*}\\
where $S(t)$ denotes the Stokes semigroup and $\mathbb{P}$ denotes the Helmholtz projection. Since $\mathbb{P}(u\cdot \nabla u)=\mathbb{P}(\omega u^{\perp})$ for $u^{\perp}={}^{t}(-u^{2},u^{1})$, $(\omega,u)$ satisfies

\begin{align*}
&\omega(t)=T(t)\omega_0+\int_{0}^{t}\nabla^{\perp}\cdot S(t-s)\mathbb{P}(\omega u^{\perp})\dd s,\quad u=K\omega,   \tag{1.9} 
\end{align*}\\
for $T(t)=-\nabla^{\perp}\cdot S(t)K$. The equations (1.9) may be viewed as an integral form of the vorticity equations (1.3). Since $S(t)u_0$ is defined for $u_0\in L^{2,\infty}_{\sigma}$, $T(t)$ is defined for all $\omega_0\in M$. We show that by the Green matrix of $S(t)$, $T(t)$ is represented by  

\begin{equation*}
\begin{aligned}
T(t)\omega_0&=\int_{\mathbb{R}^{2}_{+}}W(x,y,t)\omega_0(\dd y),\\
W(x,y,t)&=\Gamma(x-y,t)+\Gamma(x-y^{*},t)+4\int_{0}^{y_2}\int_{\mathbb{R}}\Gamma (x-z^{*}, t)\partial_{z_1}^{2}E(z-y)\dd z\\
&\quad -2\Gamma_0(x_2,t)\int_{\mathbb{R}}\Gamma_0(x_1-z_1,t)P_{y_2}(z_1-y_1)\dd z_1,\\
\Gamma_0(r,t)&=\frac{1}{(4\pi t)^{1/2}}e^{-r^{2}/4t}.
\end{aligned}
\tag{1.10}
\end{equation*}\\
With the operators,  

\begin{equation*}
\begin{aligned}
e^{t\Delta_N}\varphi&=\int_{\mathbb{R}^{2}_{+}}\left(\Gamma(x-y,t)+\Gamma(x-y^{*},t)\right)\varphi(y)\dd y,\\
e^{t\Delta_D}\varphi&=\int_{\mathbb{R}^{2}_{+}}\left(\Gamma(x-y,t)-\Gamma(x-y^{*},t)\right)\varphi(y)\dd y, \\ 
e^{t\partial_1^{2}}g&=\int_{\mathbb{R}}\Gamma_0(x_1-y_1,t)g(y_1)\dd y_1.\\
(-\Delta_D)^{-1}\omega_0&=\int_{\mathbb{R}^{2}_{+}}D(x,y)\omega_0(\dd y),
\end{aligned}
\tag{1.11}
\end{equation*}\\
$T(t)$ is represented by

\begin{equation*}
\begin{aligned}
T(t)\omega_0=
e^{t\Delta_N}\omega_0-(e^{t\Delta_N}-e^{t\Delta_D})(H\partial_2-\partial_1)\partial_1(-\Delta_D)^{-1}\omega_0-2\Gamma_0(x_2,t)e^{t\partial_1^{2}}u^{1}_{0}(\cdot,0).
\end{aligned}
\tag{1.12}
\end{equation*}\\
If (1.5) is satisfied, $T(t)\omega_0$ agrees with $T_0(t)\omega_0$ (see Theorem A.3). But the kernel $W(x,y,t)$ is different from $W_0(x,y,t)$. The formulas (1.10) and (1.12) are available to represent vorticity of the Stokes flow even if (1.5) is not satisfied. 

An important property of the operator $T(t)$ is the $L^{1}$-estimate

\begin{align*}
||T(t)\omega_0||_{1}\leq C||\omega_0||_{M},\quad t>0.
 \tag{1.13}
\end{align*}\\
This follows from integrability of the kernel $W(x,y,t)$ for the $x$-variable. We shall show that the kernel

\begin{align*}
W^{*}(x,y,t)=\Gamma(x-y^{*},t)+4\int_{0}^{y_2}\int_{\mathbb{R}}\Gamma (x-z^{*}, t)\partial_{z_1}^{2}E(z-y)\dd z
\end{align*}\\
agrees with $-G_{11}^{*}(y,x,t)$ for the Green matrix $G_{ij}(x,y,t)=\Gamma(x-y,t)\delta_{ij} +G^{*}_{ij}(x,y,t)$ of $S(t)$ (see Section 2 for the definition of $G_{ij}(x,y,t)$). Note that in contrast to $\mathbb{R}^{2}$, $S(t)$ does not satisfy the $L^{1}$-estimate \cite{DHP}, \cite{Saal07}, i.e., 

\begin{align*}
S(1)v_0\notin L^{1}(\mathbb{R}^{2}_{+})\quad \textrm{for some}\ v_0\in L^{2}\cap L^{1}(\mathbb{R}^{2}_{+}),\ \D\ v_0=0,\ v^{2}_{0}(x_1,0)=0.
\end{align*}\\
On the other hand,  since $G_{ij}(x,y,t)$ satisfies a Gaussian bound for the $y_2$-variable, the $L^{\infty}$-estimate 

\begin{align*}
||S(t)v_0||_{\infty}
\leq C||v_0||_{\infty},\quad t>0,   \tag{1.14}
\end{align*}\\
holds \cite{DHP}, \cite{Sl03}. Since $W^{*}(x,y,t)=-G_{11}^{*}(y,x,t)$, (1.13) is obtained similarly to (1.14) and is different from the $L^{1}$-estimate of $S(t)$.

The continuity at $t=0$ depends on initial conditions. We set 

\begin{align*}
C_0(\overline{\mathbb{R}^{2}_{+}})=\left\{\varphi\in C(\overline{\mathbb{R}^{2}_{+}})\ \middle|\ \lim_{|x|\to\infty}\varphi(x)=0   \right\}.
\end{align*}\\
The space $C_0$ is the pre-dual space of $M$ \cite{Rudin}. We consider the vague (weak-star) topology on $M$. Let $\delta_{0,[0,\infty)}$ denote the Dirac measure on $[0,\infty)$ at $x_2=0$, i.e., 

\begin{align*}
<\delta_{0,[0,\infty)},\psi>=\psi(0),\quad \psi\in C_0[0,\infty),
\end{align*}\\
where $<\cdot,\cdot>$ denotes the paring for $M[0,\infty)$ and $C_0[0,\infty)$. For $\omega_0\in M$, we shall show that 

\begin{align*}
T(t)\omega_0 \to \omega_0-\delta_{0,[0,\infty)}u^{1}_{0}(x_1,0)  \quad  \textrm{vaguely on}\ M\quad \textrm{as}\ t\to0.\tag{1.15}
\end{align*}\\
Since $\nabla^{\perp}_{x}D(x,y)=0$ for $y_2=0$, $K\delta_{0,[0,\infty)}u^{1}_{0}(\cdot,0)\equiv 0$. Thus by normalizing $\omega_0\in M$ by 

\begin{align*}
\tilde{\omega}_0=\omega_0-\delta_{0,[0,\infty)}u^{1}_{0}(x_1,0),  \tag{1.16}
\end{align*}\\
(1.15) is rephrased as $T(t)\tilde{\omega}_0=T(t)\omega_0\to \tilde{\omega}_0$ vaguely on $M$ as usual. Since $T(t)\omega_0$ becomes vaguely continuous by the normalization, we simply say that $T(t)\omega_0$ is vaguely continuous on $M$ at $t=0$. If $\omega_0$ is a continuous measure, $t^{1-1/q}T(t)\omega_0$ tends to zero on $L^{q}$ for $q\in (1,\infty]$. 

If $\omega_0$ has a density (i.e., $\omega_0\in L^{1}$) and (1.5) is satisfied, the stronger convergence $T(t)\omega_0\to \omega_0$ on $L^{1}$ holds. The condition (1.5) is necessary for the $L^{1}$-convergence since the trace of $S(t)u_0$ vanishes for $t>0$ and 

\begin{align*}
||S(t)u_0-u_0||_{L^{2,\infty}(\mathbb{R}^{2}_{+})}+||S(t)u_0-u_0||_{L^{1}(\partial\mathbb{R}^{2}_{+})}
+\left|\int_{\mathbb{R}^{2}_{+}}(T(t)\omega_0-\omega_0)\dd y\right|
\lesssim ||T(t)\omega_0-\omega_0||_{L^{1}(\mathbb{R}^{2}_{+})},
\end{align*}\\
by the Biot-Savart law (1.2). The $L^{1}$-convergence of $T(t)\omega_0$ also implies the zero total mass for $\omega_0$ and continuity of $S(t)u_0$ on $L^{2,\infty}$.

We construct solutions of the vorticity equations (1.3) for $u_0\in L^{2,\infty}_{\sigma}$ satisfying $\omega_0\in M$ with a small pure point part. We say that a measure $\mu\in M$ is pure point (discrete) if there exists a countable set $\{x_j\}\subset \overline{\mathbb{R}^{2}_{+}}$ and $\{\kappa_j\}\subset \mathbb{R}$ such that $\mu=\sum_{j}\kappa_j\delta_{x_j}$. A measure $\mu\in M$ is called continuous if $\mu(\{x\})=0$ for $x\in \overline{\mathbb{R}^{2}_{+}}$. If the total variation of $\mu$ is finite, the set $D=\{x\in \overline{\mathbb{R}^{2}_{+}}\ |\ \mu(\{x\})\neq 0\} $ is countable. Hence, $\mu\in M$ is uniquely decomposed as

\begin{align*}
\mu=\mu_{\textrm{pp}}+\mu_{\textrm{cont}}
\end{align*}\\
with pure point $\mu_{\textrm{pp}}$ and continuous $\mu_{\textrm{cont}}$ by setting $\mu_{\textrm{pp}}(E)=\mu(D\cap E)$ for Borel sets $E\subset \overline{\mathbb{R}^{2}_{+}}$. Since $t^{1-1/q}T(t)\omega_{0,\textrm{pp}}$ does not tend to zero as $t\to 0$, we assume a smallness for $\omega_{0,\textrm{pp}}$ in order to construct local-in-time solutions. If the total variation of $\omega_0\in M$ is small, we are able to construct small global-in-time solutions. Let $BC([0,T]; X)$ (resp. $BC_{w}([0,T]; X)$) denote the space of all bounded (resp. weakly-star) continuous functions from $[0,T]$ to a Banach space $X$. We denote by $BC((0,T]; X)$ the space of all bounded functions in $[0,T]$, continuous in $(0,T]$. The main result of this paper is:     

\vspace{5pt}

%thm1.1
\begin{thm}
(i) There exists $\delta_1>0$ such that for $u_0\in L^{2,\infty}_{\sigma}$ satisfying $\omega_0\in M$ and $||\omega_{0,\textrm{pp}}||_{M}\leq \delta_1$, there exists $T>0$ and a unique $(\omega,u)$ satisfying (1.8), (1.9) and

\begin{align*}
&\omega\in BC_{w}([0,T]; M),\tag{1.17} \\
&u\in BC_{w}([0,T]; L^{2,\infty}),\tag{1.18} \\
&t^{1-1/q}\omega\in BC((0,T]; L^{q}), \quad 1<q\leq  \infty,  \tag{1.19}\\  
&t^{1/2-1/p}u\in BC((0,T]; L^{p}),\quad 2<p\leq \infty.  \tag{1.20}
\end{align*}\\
If $\omega_{0}$ is continuous, both values (1.19) and (1.20) vanish at $t=0$. If in addition that $\omega_0\in L^{1}$ and $u^{1}_{0}(x_1,0)=0$, $(\omega,u)$ is strongly continuous at $t=0$.

\noindent
(ii) There exists $\delta_2>0$ such that for $u_0\in L^{2,\infty}_{\sigma}$ satisfying $\omega_0\in M$ and $||\omega_0||_{M}\leq \delta_2$, there exists a unique $(\omega,u)$ satisfying (1.8), (1.9), (1.17)-(1.20) for $T=\infty$. 
\end{thm}

\vspace{5pt}

Since (1.1) is globally well-posed for bounded initial data with finite Dirichlet integral \cite{A7}, by replacing $t\in (0,T]$ as an initial time, we have:

\vspace{5pt}

%thm1.2
\begin{thm}
The solution constructed in Theorem 1.1 (i) is global, i.e., $(\omega,u)$ satisfies (1.8), (1.9), (1.17)-(1.20) for all $T>0$.
\end{thm}

\vspace{15pt}

Theorem 1.2 implies global well-posedness of (1.1) for $\omega_0\in M$ with a small pure point part (e.g., $\omega_{0,\textrm{pp}}\equiv 0$). It in particular implies that vortex sheets diffuse by the Navier-Stokes flow with boundary. On the other hand, smallness conditions are assumed in Theorems 1.1 (ii) and 1.2  for the pure point part $\omega_{0,\textrm{pp}}$ in order to construct global-in-time solutions. Existence for $\omega_0\in M$ with large $\omega_{0,\textrm{pp}}$ is unknown even if $\omega_0$ is a point mass, i.e., $\omega_0=\kappa \delta_{x_0}$ for $x_0\in \mathbb{R}^{2}_{+}$ and large $\kappa \in \mathbb{R}$. For the Stokes flow, $\kappa W(x,x_0,t)$ defined by $(1.10)_2$ is an exact solution for $\omega_0=\kappa\delta_{x_0}$.

For the Cauchy problem, global-in-time solutions of (1.1) exist for all $\omega_0\in M(\mathbb{R}^{2})$ by a priori estimates of vorticity \cite{Cottet}, \cite{GMO}. The uniqueness for $\omega_0\in M(\mathbb{R}^{2})$ with small $\omega_{0,\textrm{pp}}$ is proved in \cite{GMO} based on an integral form of the vorticity equations. See also \cite{Kato94}. The uniqueness for $\omega_0\in M(\mathbb{R}^{2})$ with large $\omega_{0,\textrm{pp}}$ is more difficult. For $u_0=(2\pi)^{-1}x^{\perp}|x|^{-2}\in L^{2,\infty}(\mathbb{R}^{2})$ with $\omega_0=\delta_0\in M(\mathbb{R}^{2})$, there exists a forward self-similar solution of (1.1) in $\mathbb{R}^{2}$, called the Lamb-Oseen vortex:

\begin{align*}
\Omega(x,t)=\frac{1}{t}\Omega_0\left(\frac{x}{\sqrt{t}}\right),
\quad U(x,t)=\frac{1}{\sqrt{t}}U_0\left(\frac{x}{\sqrt{t}}\right), 
\end{align*}\\
where

\begin{align*}
\Omega_0(x)=\frac{1}{4\pi}e^{-|x|^{2}/4},
\quad U_0(x)=\frac{x^{\perp}}{2\pi |x|^{2}}\left(1-e^{-|x|^{2}/4}\right). 
\end{align*}\\
The uniqueness for $\omega_0=\kappa \delta_0$ and large $\kappa\in \mathbb{R}$ is proved in \cite{GallayWayne} by using a relative entropy for the self-similar transform of $\omega$. See also \cite{GGL} for an alternative proof. The uniqueness relates to the asymptotic formula 

\begin{align*}
\lim_{t\to \infty}t^{1-1/q}||\omega-\kappa \Omega||_{q}=0,\quad \omega_0\in M(\mathbb{R}^{2}),\  \kappa =\int_{\mathbb{R}^{2}}\omega_0(\dd y),\ 1\leq q\leq \infty.  \tag{1.21}
\end{align*}\\
The formula (1.21) is studied in \cite{KG} for $\omega_0$ with a small total variation and extended in \cite{Carpio} for small $\kappa$ (see also \cite{GW02}, \cite{GGS}). For large $\kappa$, (1.21) is proved in \cite{GallayWayne}. The uniqueness for general $\omega_0\in M(\mathbb{R}^{2})$ with large $\omega_{0,\textrm{pp}}$ is proved in \cite{GG05}. 

For the half plane, initial data $u_0$ of homogeneous of degree $-1$ satisfying $\D\ u_0=0$ and $u^{2}_0(x_1,0)=0$ are only of the form 

\begin{align*}
u_0(x)=\kappa(\theta)\frac{x}{|x|^{2}},
\end{align*}\\
for some $\kappa(\theta)$, due to the boundary condition. Here, $(r,\theta)$ is the polar coordinate. Obviously, $\omega_0(x)=\kappa(\theta)|x|^{-2}\notin M$. Hence, any forward self-similar solutions of (1.1) in $\mathbb{R}^{2}_{+}$ do not satisfy the initial condition $\omega_0\in M$, in contrast to $\mathbb{R}^{2}$. As noted in \cite{GallayWayne}, there are forward self-similar solutions in $\mathbb{R}^{2}$ such that $\omega_0\notin M$. For the half plane, existence of small forward self-similar solutions follows from a result in $L^{2,\infty}$ \cite{KY}. 

It is an interesting question whether solutions for $u_0\in L^{2,\infty}_{\sigma}$ satisfying $\omega_0\in M$ tend to zero as time goes to infinity. For the Stokes flow, we have

\begin{align*}
\lim_{t\to\infty} t^{1-1/q}||T(t)\omega_0||_{q}=0,\quad \omega_0\in M,\ 1\leq q\leq \infty.  \tag{1.22}
\end{align*}\\
See Theorem 4.7. If the total variation of $\omega_0$ is small, $t^{1-1/q}\omega$ is globally bounded in $L^{q}$ by Theorem 1.1 (ii). It is unknown whether $t^{1-1/q}\omega$ tend to zero as $t\to\infty$. If $u_0\in L^{2}$, we have $\lim_{t\to\infty}||u||_{2}=0$ \cite{BM88}. The large time behavior is important to study non-existence of backward solutions. We refer to \cite{Seregin15} for a Liouville theorem in $L^{\infty}(-\infty,0; L^{2})$.

\vspace{10pt}
This paper is organized as follows. In Section 2, we prove the Biot-Savart law (1.2). In Section 3, we prove the formulas (1.10), (1.12) and a kernel estimate for $W(x,y,t)$. In Section 4, we study continuity of $T(t)$ at time zero. We also prove the asymptotic formula (1.22). In Section 5, we prove Theorems 1.1 and 1.2. In Appendix A, we give a proof for the formula (1.4).

\vspace{15pt}

%section2
\section{The Stokes flow on $L^{2,\infty}$}

\vspace{15pt}

In this section, we prove the Biot-Savart law (1.2) for solenoidal vector fields $u_0\in L^{2,\infty}(\mathbb{R}^{2}_{+})$ with a finite Borel measure $\omega_0=-\nabla^{\perp}\cdot u_0$ on $\overline{\mathbb{R}^{2}_{+}}$ (Lemma 2.3). We define all function spaces used in the subsequent sections.

\vspace{15pt}

\subsection{Solenoidals in $L^{2,\infty}$}
We recall the Lorentz space \cite{SteinWeiss}, \cite{BL}, \cite{Ad}. For a measurable function $f$ in $\mathbb{R}^{2}_{+}$ we set a distribution function $m(t,f)$ and a decreasing rearrangement $f^{*}(t)$ by 

\begin{align*}
&m(t,f)=\left|\left\{x\in \mathbb{R}^{2}_{+}\ |\ |f(x)|>t \right\}\right|, \\
&f^{*}(t)=\inf \{s\in (0,\infty)\ |\ m(s)<t\},\quad t>0,
\end{align*}\\ 
where $|E|$ denotes the Lebesgue measure for a measurable set $E\subset \mathbb{R}^{2}_{+}$. For $p\in (1,\infty)$, we define $L^{p,q}(\mathbb{R}^{2}_{+})$ by the space of all measurable functions $f$ such that 

\begin{align*}
&||f||_{p,q}=\left(\int_{0}^{\infty}(t^{1/p}f^{*}(t))^{q}\frac{\dd t}{t}\right)^{1/q}<\infty,\quad 1\leq q<\infty,\\
&||f||_{p,\infty}=\sup_{t>0}t^{1/p}f^{*}(t) <\infty,\quad q=\infty.
\end{align*}\\
The space $L^{p,q}$ agrees with $L^{p}$ if $q=p$ and $L^{p,q_1}\subset L^{p,q_2}$ for $1\leq q_1\leq q_2\leq \infty$ \cite[p.16]{BL}. In particular, $L^{p}\subset L^{p,\infty}$. A function $f$ belongs to $L^{p,\infty}$ if and only if 

\begin{align*}
\sup_{E}|E|^{-1+1/p}\int_{E}|f(x)|\dd x<\infty.
\end{align*}\\
This becomes an equivalent norm to $||\cdot||_{p,\infty}$ \cite{GM89}. The space $L^{p,q}$ is a quasi-normed Banach space and agrees with the real interpolation space $(L^{p_0}, L^{p_1})_{\theta,q}$, i.e., 

\begin{align*}
(L^{p_0},L^{p_1})_{\theta, q}=L^{p,q},\quad 1< p_0<p< p_1< \infty,\ 1/p=(1-\theta)/p_0+\theta/p_1.
\end{align*}\\
By a duality theorem \cite[3.7.1 Theorem]{BL}, 

\begin{align*}
(L^{p,q})^{*}=(L^{p_0},L^{p_1})_{\theta,q}^{*}
=(L^{p_0'},L^{p_1'})_{\theta,q'}
=L^{p',q'},\quad 1\leq q<\infty,
\end{align*}\\
where $p'$ denotes the conjugate exponent to $p$. We denote by $C_{c}^{\infty}(\mathbb{R}^{2}_{+})$ the space of all smooth functions with compact support in $\mathbb{R}^{2}_{+}$. Since $C_{c}^{\infty}$ is dense in $L^{p_0}\cap L^{p_1}$, $C^{\infty}_{c}$ is also dense in $L^{p,q}$ for $1\leq q<\infty$ \cite[3.4.2 Theorem (b)]{BL}. In the sequel, we do not distinguish the space of scaler and vector-valued functions.

We set the subspaces of $L^{p}$ by

\begin{align*}
&L^{p}_{\sigma}(\mathbb{R}^{2}_{+})=\overline{C_{c,\sigma}^{\infty}(\mathbb{R}^{2}_{+})}^{||\cdot||_{L^{p}}},\ C_{c,\sigma}^{\infty}(\mathbb{R}^{2}_{+})=\left\{f\in C_{c}^{\infty}(\mathbb{R}^{2}_{+})\ \middle|\ \D\ f=0 \right\},\\
&G^{p}(\mathbb{R}^{2}_{+})=\left\{\nabla \Phi\in L^{p}(\mathbb{R}^{2}_{+})\ \middle|\ \Phi\in L^{1}_{\textrm{loc}}(\mathbb{R}^{2}_{+}) \right\}.
\end{align*}\\
The space $L^{p}$ is decomposed into the direct sum 

\begin{align*}
L^{p}(\mathbb{R}^{2}_{+})=L^{p}_{\sigma}(\mathbb{R}^{2}_{+})
\oplus G^{p}(\mathbb{R}^{2}_{+}).
\end{align*}\\
We call $\mathbb{P}:L^{p} \longrightarrow L^{p}_{\sigma}$ the Helmholtz projection operator (e.g., \cite{BM88}). The space $L^{p}_{\sigma}$ agrees with the space of all $L^{p}$-solenoidal vector fields in $\mathbb{R}^{2}_{+}$, i.e., 

\begin{align*}
L^{p}_{\sigma}(\mathbb{R}^{2}_{+})
=\left\{ f\in L^{p}(\mathbb{R}^{2}_{+})\ \middle|\ \D\ f=0\ \textrm{in}\ \mathbb{R}^{2}_{+},\ f^{2}(x_1,0)=0,\ x_1\in \mathbb{R}\  \right\}.  \tag{2.1}
\end{align*}\\
The normal trace $f^{2}(x_1,0)$ is understood in the Sobolev space of a negative order $W^{-1/p',p'}(\mathbb{R})=W^{1-1/p,p}(\mathbb{R})^{*}$ \cite[II. 1.2.3 Lemma]{Sohr}, \cite[Theorem II 10.2]{Gal}. Indeed, for $f\in L^{p}$ satisfying $\D\ f=0$ and $f^{2}(x_1,0)=0$, set 

\begin{align*}
f=\mathbb{P}f+(I-\mathbb{P})f=f_0+\nabla \Phi.
\end{align*}\\
Since $\Delta \Phi=0$ and $\partial_2 \Phi(x_1,0)=0$, $\nabla \Phi\equiv 0$ and $f=\mathbb{P}f\in L^{p}_{\sigma}$ by the Liouville theorem. This implies (2.1).

Following \cite{MY92}, \cite{KY}, \cite{Yamazaki}, we define the $L^{p,q}$-solenoidal space. For the two interpolation pairs $\{L^{p_0},L^{p_1}\}$ and $\{L^{p_0}_{\sigma},L^{p_1}_{\sigma}\}$, $\mathbb{P}: L^{p_i}\longrightarrow L^{p_i}_{\sigma}$, $i=0,1$, is bounded and surjective. Since $(L^{p_0}, L^{p_1})_{\theta,q}$ and $(L^{p_0}_{\sigma}, L^{p_1}_{\sigma})_{\theta,q}$ are exact interpolation spaces of type $\theta$ \cite[7.23 THEOREM]{Ad},

\begin{align*}
\mathbb{P}: L^{p,q}=(L^{p_0},L^{p_1})_{\theta,q}\longrightarrow (L^{p_0}_{\sigma},L^{p_1}_{\sigma})_{\theta,q}\subset L^{p,q}
\end{align*}\\
is bounded and surjective. We set $L^{p,q}_{\sigma}:=\mathbb{P}L^{p,q}=(L^{p_0}_{\sigma},L^{p_1}_{\sigma})_{\theta,q}$. Since $C_{c,\sigma}^{\infty}$ is dense in $L^{p_0}_{\sigma}\cap L^{p_1}_{\sigma}$, $C_{c,\sigma}^{\infty}$ is dense in $L^{p,q}_{\sigma}$ for $1\leq q<\infty$. Moreover, we have

\begin{align*}
L^{p,q}_{\sigma}(\mathbb{R}^{2}_{+})=\left\{f\in L^{p,q}(\mathbb{R}^{2}_{+})\ \middle|\ \D\ f=0\ \textrm{in}\ \mathbb{R}^{2}_{+},\ f^{2}(x_1,0)=0,\ x_1\in \mathbb{R} \right\}.  \tag{2.2}
\end{align*}\\
Since $L^{p,q}_{\sigma}= (L^{p_0}_{\sigma},L^{p_1}_{\sigma})_{\theta,q}\subset L^{p_0}_{\sigma}+L^{p_1}_{\sigma}$ by the definition of the real interpolation, the right-hand side of (2.2) is larger than the left-hand side. The converse inclusion follows in the same way as (2.1).

For $p=\infty$, we define $L^{\infty}_{\sigma}$ by the space of all $f\in L^{\infty}$ satisfying $\D\ f=0$ and $f^{2}(x_1,0)=0$ \cite{DHP}.

\vspace{15pt}

\subsection{The Biot-Savart law}

We recall the space of finite real regular Borel measures on $\overline{\mathbb{R}^{2}_{+}}$ \cite{Rudin}, \cite{Fo}. Let ${\mathcal{B}}$ be a Borel $\sigma$-algebra on $\overline{\mathbb{R}^{2}_{+}}$ (i.e., the $\sigma$-algebra generated by open sets in $\overline{\mathbb{R}^{2}_{+}}$). We say that $\mu: \mathcal{B}\longrightarrow [-\infty,\infty]$ is a singed Borel measure if $\mu$ is countably additive. For positive $\mu\geq 0$, $E\in {\mathcal{B}}$ is called outer (resp. inner) regular if

\begin{align*}
\mu(E)=\inf\ \left\{\mu(U)\ |\ E\subset U,\ U:\textrm{open}  \right\}
\end{align*}\\
(resp. $\mu(E)=\sup\ \{\mu(K)\ |\ K\subset E,\ K:\textrm{compact}  \}$). If all $E\in \mathcal{B}$ are outer and inner regular, $\mu$ is called regular. By the Jordan decomposition $\mu=\mu_{+}-\mu_{-}$ for $\mu_+,\mu_{-}\geq 0$, we set the total variation measure 

\begin{align*}
|\mu|=\mu_{+}+\mu_{-}.
\end{align*}\\
A signed measure $\mu$ is called regular if $|\mu|\geq 0$ is regular. We denote by $M(\overline{\mathbb{R}^{2}_{+}})$ the space of all signed regular Borel measures on $\overline{\mathbb{R}^{2}_{+}}$ equipped with the norm $||\mu ||_{M}=|\mu|(\overline{\mathbb{R}^{2}_{+}})$. We set 

\begin{align*}
C_{0}(\overline{\mathbb{R}^{2}_{+}})
=\left\{\varphi\in C(\overline{\mathbb{R}^{2}_{+}})\ \middle|\ \lim_{|x|\to\infty}\varphi(x)=0   \right\}.
\end{align*}\\
The space $C_c^{\infty}(\overline{\mathbb{R}^{2}_{+}})$ is dense in $C_{0}(\overline{\mathbb{R}^{2}_{+}})$. By the Riesz representation theorem, 

\begin{align*}
M(\overline{\mathbb{R}^{2}_{+}})=C_{0}(\overline{\mathbb{R}^{2}_{+}})^{*}.
\end{align*}\\
The weak-star topology of $M(\overline{\mathbb{R}^{2}_{+}})$ is called vague topology \cite{Fo}. In the sequel, we write $M=M(\overline{\mathbb{R}^{2}_{+}})$ by omitting the symbol $\overline{\mathbb{R}^{2}_{+}}$.

\vspace{15pt}

%prop2.1
\begin{prop}
Set 

\begin{align*}
K\mu (x)=\int_{\mathbb{R}^{2}_{+}}\nabla^{\perp}_{x}D(x,y)\mu(\dd y).  \tag{2.3}
\end{align*}\\
\noindent 
Then, 
\begin{align*}
K: 
\begin{cases}
&M \longrightarrow L^{2,\infty}, \\
&L^{q}\longrightarrow L^{p},\ q\in (1,2),\ 1/p=1/q-1/2,\\
&L^{2,1}\longrightarrow C_{0},
\end{cases}
\tag{2.4}
\end{align*}
is bounded. 
\end{prop}

\vspace{5pt}

\begin{proof}
For a measure $\mu\in M$, we set a measure $\bar{\mu}$ on $\mathbb{R}^{2}$ by

\begin{align*}
\bar{\mu}(E)=\mu(E\cap \overline{\mathbb{R}^{2}_{+}})
\end{align*}\\
for Borel sets $E\subset \mathbb{R}^{2}$. By the reflection $\bar{\mu}^{*}(E)=\bar{\mu}(E^{*})$ and $E^{*}=\{x^{*} |  x\in E\}$, we set $\mu^{\textrm{odd}}=\bar{\mu}-\bar{\mu}^{*}$. By changing the variable, 

\begin{align*}
K\mu(x)
=\int_{\mathbb{R}^{2}_{+}}(\nabla^{\perp}_{x}E(x-y)-\nabla^{\perp}_{x}E(x-y^{*}))\mu(\dd y)
=\int_{\mathbb{R}^{2}}\nabla^{\perp}_{x}E(x-y)\mu^{\textrm{odd}}(\dd y)=\nabla^{\perp} E*\mu^{\textrm{odd}}.
\end{align*}\\
Since $\mu^{\textrm{odd}}\longmapsto \nabla^{\perp}E*\mu^{\textrm{odd}}$ is bounded from $M(\mathbb{R}^{2})$ to $L^{2,\infty}(\mathbb{R}^{2})$ \cite[Lemma 2.2 (i)]{GMO}, $(2.4)_{1}$ follows.

For $\mu\in L^{q}$, $q\in (1,2)$, we set the odd extension  

\begin{align*}
\mu^{\textrm{odd}}(x)=
\begin{cases}
\quad \mu(x_1,x_2),& \quad x_2\geq 0,\\
\quad-\mu(x_1,-x_2),& \quad x_2<0.
\end{cases}
\tag{2.5}
\end{align*}\\
Since $\mu^{\textrm{odd}}\longmapsto \nabla^{\perp}E*\mu^{\textrm{odd}}$ is bounded from $L^{q}(\mathbb{R}^{2})$ to $L^{p}(\mathbb{R}^{2})$ \cite{GMO}, $(2.4)_2$ follows. Since $C_0^{*}=M$ and $(L^{2,1})^{*}=L^{2,\infty}$, $\mu^{\textrm{odd}}\longmapsto \nabla^{\perp}E*\mu^{\textrm{odd}}$ is bounded from $L^{2,1}(\mathbb{R}^{2})$ to $C_{0}(\mathbb{R}^{2})$. Thus $(2.4)_3$ follows.
\end{proof}

\vspace{15pt}

%prop2.2
\begin{prop}
Set 

\begin{align*}
T_1\mu(x_1)=\frac{1}{\pi}\int_{\mathbb{R}}\frac{y_2}{|x_1-y_1|^{2}+y_2^{2}}\mu(\dd y).  \tag{2.6}
\end{align*}\\
Then, 
\begin{align*}
T_1: 
\begin{cases}
&M(\overline{\mathbb{R}^{2}_{+}})\longrightarrow L^{1}(\mathbb{R}), \\
&L^{q}(\mathbb{R}^{2}_{+}) \longrightarrow L^{s,q}(\mathbb{R}),\ 1< q< 2,\ 1/(2s)=1/q-1/2,\\
&L^{2,1}(\mathbb{R}^{2}_{+}) \longrightarrow C_0(\mathbb{R}),
\end{cases}
\tag{2.7}
\end{align*}
is bounded. 
\end{prop}

\vspace{5pt}

\begin{proof}
By integrating (2.6) directly, $(2.7)_1$ follows. Since $(K\mu)^{1}(x_1,0)=T_1\mu(x_1)$, $(2.7)_3$ follows from $(2.4)_{3}$. By applying the general Marcinkiewicz interpolation theorem \cite[5.3.2 Theorem]{BL} to $(2.7)_1$ and $(2.7)_3$, $(2.7)_{2}$ follows.
\end{proof}

\vspace{15pt}

%lem2.3
\begin{lem}
(i) For $u\in L^{2,\infty}_{\sigma}$ (resp. $u\in L^{p}_{\sigma}$) satisfying $\omega\in M$ (resp. $\omega\in L^{q}$, $q\in (1,2)$, $1/p=1/q-1/2$), $u=K\omega$ holds. 

\noindent
(ii) For $\omega\in M$ and $u=K\omega$, set 

\begin{align*}
\tilde{\omega}=\omega-\delta_{0,[0,\infty)}u^{1}(x_1,0)    \tag{2.8}
\end{align*}\\
by the Dirac measure $\delta_{0,[0,\infty)}$ on $[0,\infty)$ at $x_2=0$. Then, $K\tilde{\omega}=K\omega$.
\end{lem}

\vspace{5pt}

\begin{proof}
Since $u\in L^{2,\infty}_{\sigma}$ satisfies $\nabla\cdot u=0$ and $u^{2}(x_1,0)=0$, there exists a stream function $\psi$ such that $u=\nabla^{\perp}\psi$. Since $0=u^{2}(x_1,0)=-\partial_1\psi(x_1,0)$, we may assume that $\psi(x_1,0)=0$. We set $\tilde{u}=K\omega=\nabla \tilde{\psi}\in L^{2,\infty}$ (resp. $\tilde{u}\in L^{p}$). Since $-\Delta \psi=\omega$ and $\psi(x_1,0)=0$, $\varphi=\psi-\tilde{\psi}$ satisfies $-\Delta\varphi=0$ and $\varphi(x_1,0)=0$. Since $\partial_1\varphi\in L^{2,\infty}$ (resp. $\partial_1\varphi\in L^{p}$), applying the Liouville theorem implies $\partial_1\varphi\equiv 0$. Hence, $\partial_2^{2}\varphi\equiv 0$. Since $\nabla \varphi\to 0$ as $|x|\to\infty$ and $\varphi(x_1,0)=0$, $\varphi\equiv 0$ follows. Thus, $u=K\omega$. This proves (i).

We prove (ii). By $(2.7)_1$, $u^{1}(y_1,0)\in L^{1}(\mathbb{R})$. For fixed $x\in \mathbb{R}^{2}_{+}$ and $y_1\in \mathbb{R}$ ($y_1\neq x_1$), observe that

\begin{align*}
\nabla^{\perp}_{x}D(x,y)
=\frac{1}{2\pi}\left(\frac{(x-y)^{\perp}}{|x_1-y_1|^{2}+|x_2-y_2|^{2}} 
-\frac{(x-y^{*})^{\perp}}{|x_1-y_1|^{2}+|x_2+y_2|^{2}}
\right)\in C_0[0,\infty)
\end{align*}\\
as a function of $y_2\in[0,\infty)$. Let $<\cdot,\cdot>$ denote the pairing for $M[0,\infty)$ and $C_0[0,\infty)$. Since $\nabla_{x}^{\perp}D(x,y)=0$ for $y_2=0$, it follows that 

\begin{align*}
K\delta_{0,[0,\infty)}u^{1}(\cdot,0)
=\int_{\mathbb{R}}<\delta_{0,[0,\infty)},\nabla^{\perp}D(x,y) >u^{1}(y_1,0)\dd y_1=0.
\end{align*}\\
We proved $K\tilde{\omega}=K\omega$.
\end{proof}

\vspace{15pt}

\subsection{The Stokes semigroup}

We define the Stokes semigroup 

\begin{align*}
S(t)u_0=\int_{\mathbb{R}^{2}_{+}}G(x,z,t)u_0(z)\dd z,  \tag{2.9}
\end{align*}\\
by the Green matrix $G=(G_{ij})_{1\leq i,j\leq 2}$, 

\begin{equation*}
\begin{aligned}
&G_{ij}(x,z,t)=\Gamma(x-z,t)\delta_{ij}+G_{ij}^{*}(x,z,t),\\
&G_{ij}^{*}(x,z,t)=-\Gamma(x-z^{*},t)\delta_{ij}-4(1-\delta_{j2})\partial_{x_j}\int_{0}^{x_2}\int_{\mathbb{R}}\partial_{x_i} E(x-w)\Gamma(w-z^{*},t)\dd w,
\end{aligned}
\tag{2.10}
\end{equation*}\\
and the Kronecker delta $\delta_{ij}$ \cite[p.336]{Sl03b}. The function $G^{*}$ satisfies the pointwise estimate 

\begin{align*}
|\partial_t^{s}\partial_x^{k}\partial_z^{m}G^{*}(x,z,t)|
\leq \frac{Ce^{-cz_2^{2}/t}}{t^{s+m_2/2}(x_2^{2}+t)^{k_2/2}(|x-z^{*}|^{2}+t)^{(k_1+m_1+2)/2}},  \tag{2.11}
\end{align*}\\
for $k=(k_1,k_2)$, $m=(m_1,m_2)$ and $s\geq 0$ with some constant $C$ \cite[Proposition 2.5]{Sl03b}. By (2.11), $S(t)$ satisfies the $L^{p}-L^{r}$-estimate

\begin{align*}
||\partial_t^{s}\partial_x^{k}S(t)u_0||_{r}\leq \frac{C}{t^{|k|/2+s+1/p-1/r}}||u_0||_{p},\quad  t>0,\ u_0\in L^{p}_{\sigma},\ 1<p\leq r<\infty. \tag{2.12}
\end{align*}\\
See \cite[Proposition 4.1]{BM88} for $p\in (1,\infty)$ and \cite{DHP}, \cite{Sl03} for $p=\infty$. The estimate (2.12) implies that $S(t)$ is a bounded analytic semigroup on $L^{p}_{\sigma}$. Since $L^{p,q}_{\sigma}=(L^{p_0}_{\sigma},L^{p_1}_{\sigma})_{\theta,q}$ for $1<p_0<p<p_1<\infty$, by applying an interpolation theorem \cite[7.23 Theorem]{Ad}, $S(t)$ is also a bounded analytic semigroup on $L^{p,q}_{\sigma}$ for $1\leq q\leq \infty$. Since $C_{c,\sigma}^{\infty} $ is dense in $L^{p,q}_{\sigma}$ for $1\leq q<\infty$, $S(t)$ is a $C_0$-semigroup on $L^{p,q}_{\sigma}$. By the duality $(L^{p',1})^{*}=L^{p,\infty}$, $S(t)u_0$ is weakly-star continuous on $L^{p,\infty}$ at $t=0$. Moreover, we have

\begin{align*}
||\partial_t^{s}\partial_x^{k}S(t)u_0||_{r}\leq \frac{C}{t^{|k|/2+s+1/p-1/r}}||u_0||_{p,\infty},\quad t>0,\  u_0\in L^{p,\infty}_{\sigma},\ p< r\leq \infty.
\tag{2.13}
\end{align*}\\
The estimate (2.13) follows from (2.12) by taking $1<p_0<p<p_1<r$ and applying an interpolation theorem \cite[7.23 Theorem]{Ad} for $\partial_t^{s}\partial_x^{k} S(t): L^{p_i}_{\sigma}\longrightarrow L^{r}$, $i=0,1$.

To study solutions of (1.8) and (1.9), we use composition operators. The estimate (2.11) yields

\begin{align*}
||\partial S(t)\mathbb{P}f||_{p}\leq \frac{C}{t^{1/2}}||f||_{p},\quad t>0,\ f\in L^{p},\ 1<p<\infty.  \tag{2.14}
\end{align*}\\
The adjoint operator satisfies

\begin{align*}
||S(t)\mathbb{P}\partial F||_{p}\leq \frac{C}{t^{1/2}}||F||_{p},\quad t>0,\ F\in L^{p},\ 1<p<\infty, \tag{2.15}
\end{align*}\\
where $\partial=\partial_{x}^{k}$ indiscriminately denotes the spatial derivatives $|k|=1$. The operators $\partial S(t)\mathbb{P}$ and $S(t)\mathbb{P}\partial$ are understood as one operators acting on $L^{p}$. It still acts as a bounded operator for $p=1$ and $p=\infty$ even if $\mathbb{P}$ is unbounded.

\vspace{15pt}

%prop2.4
\begin{lem}
The operators $\partial S(t)\mathbb{P}$ and $S(t)\mathbb{P}\partial$ are uniquely extendable to bounded operators on $L^{1}$ and $C_{0}$ together with 

\begin{align*}
&||\partial S(t)\mathbb{P}f||_{1}
\leq \frac{C}{t^{1/2}}||f||_{1},\quad t>0,\ f\in L^{1},  \tag{2.16} \\
&||S(t)\mathbb{P}\partial F||_{\infty}
\leq \frac{C}{t^{1/2}}||F||_{\infty},\quad t>0,\ F\in C_0. \tag{2.17}  \end{align*}
\end{lem}

\vspace{5pt}

\begin{proof}
Let $(\cdot,\cdot)$ denote the pairing for $L^{1}$ and $C_0$. By integration by parts, observe that 

\begin{align*}
(\partial S(t)\mathbb{P}f, F )
=-(f, S(t)\mathbb{P}\partial F ),
\quad f,F\in C^{\infty}_{c}(\overline{\mathbb{R}^{2}_{+}}).  \tag{2.18}
\end{align*}\\
Since (2.17) holds for $F\in C^{\infty}_{c}(\mathbb{R}^{2}_{+})$ \cite{Sl03}, we estimate 

\begin{align*}
\left|(\partial S(t)\mathbb{P}f, F )\right|
=\left|(f, S(t)\mathbb{P}\partial F )\right|
\leq \frac{C}{t^{1/2}}||f||_{L^{1}(\mathbb{R}^{2}_{+})} ||F||_{L^{\infty}(\mathbb{R}^{2}_{+})}.
\end{align*}\\
By taking a supremum for $F\in C^{\infty}_{c}(\mathbb{R}^{2}_{+})$, we obtain (2.16) for $f\in C^{\infty}_{c}(\mathbb{R}^{2}_{+})$. By taking the closure in $L^{1}$, $\partial S(t)\mathbb{P}$ is uniquely extendable to a bounded operator on $L^{1}$ together with (2.16) and (2.18). By (2.18), (2.17) holds for $F\in C^{\infty}_{c}(\overline{\mathbb{R}^{2}_{+}})$. By taking the closure in $C_0$, $S(t)\mathbb{P}\partial$ is uniquely extendable to $C_0$ together with (2.17). 
\end{proof}

\vspace{5pt}

\begin{rem}
The estimate 

\begin{align*}
||\partial S(t)f||_{{\mathcal{H}}^{1}(\mathbb{R}^{2}_{+})}\leq \frac{C}{t^{1/2}}||f||_{L^{1}(\mathbb{R}^{2}_{+})}, \quad t>0,\ f\in L^{1},\ \D\ f=0,\ f^{2}(x_1,0)=0,  \tag{2.19}
\end{align*}\\
is known to hold \cite{GMSim} for the Hardy space ${\mathcal{H}}^{1}(\mathbb{R}^{2}_{+})$. The estimate (2.16) holds even if $f$ is not solenoidal. On the other hand, (2.16) is weaker than (2.19) since ${\mathcal{H}}^{1}(\mathbb{R}^{2}_{+})\subset L^{1}(\mathbb{R}^{2}_{+})$.  
\end{rem}

\vspace{15pt}

%Section 4
\section{Vorticity associated with the Stokes flow}

\vspace{15pt}

We first derive (1.12) by calculating a kernel of $T(t)=-\nabla^{\perp}\cdot S(t)K$. The explicit form of the kernel (1.10) follows from a computation of the kernel of $(H\partial_2-\partial_1)\partial_1(-\Delta_D)^{-1}$. A key fact is that the Hilbert transform of the Poisson kernel $P_s$ is the conjugate Poisson kernel $Q_{s}$, i.e., $HP_s=Q_s$, $s>0$. By using this fact, we calculate the discontinuous kernel $H\partial_2 E$. 

\vspace{15pt}

\subsection{The Hilbert transform}

To prove (1.12), we use the Hilbert transform \cite[Chapter III-VI]{SteinWeiss}. For a rapidly decreasing function $\varphi$, we set the Fourier transform 

\begin{align*}
\hat{\varphi}(\xi_1)={\mathcal{F}}\hspace{-1pt}\varphi(\xi_1)=\int_{\mathbb{R}}e^{-ix_1\xi_1}\varphi(x_1)\dd x_1, \quad i=\sqrt{-1}.
\end{align*}\\
For a tempered distribution $\varphi$, the Fourier transform is defined by 

\begin{align*}
(\hat{\varphi},\psi)=(\varphi,\hat{\psi})
\end{align*}\\
by the pairing $(\cdot,\cdot,)$ and rapidly decreasing functions $\psi$. We define the Hilbert transform $H$ by 

\begin{align*}
\widehat{H\varphi}(\xi_1)=-i\frac{\xi_1}{|\xi_1|}\hat{\varphi}(\xi_1).  \tag{3.1}
\end{align*}\\
The operator $H$ satisfies $H^{2}=-I$. It acts as a bounded operator on $L^{q}(\mathbb{R})$ for $q\in (1,\infty)$ \cite{SteinWeiss}. We set the Poisson kernel $P_s$ and the conjugate Poisson kernel $Q_s$ by

\begin{align*}
P_{s}(x_1)=\frac{s}{\pi(x_1^{2}+s^{2})},\quad Q_{s}(x_1)=\frac{x_1}{\pi(x_1^{2}+s^{2})}.  
\end{align*}\\
Their Fourier transforms are

\begin{align*}
\hat{P}_s(\xi_1)=e^{-s|\xi_1|},\quad \hat{Q}_s(\xi_1)=-i\frac{\xi_1}{|\xi_1|}e^{-s|\xi_1|}.
\end{align*}\\
We set the Poisson semigroup by $\widehat{e^{sA}\varphi}=e^{-s|\xi_1|}\hat{\varphi}$, i.e., 

\begin{align*}
e^{sA}\varphi(x_1)=(P_{s}*_1\varphi)(x_1)=\frac{1}{\pi}\int_{\mathbb{R}}\frac{s}{(x_1-y_1)^{2}+s^{2}}\varphi(y_1)\dd y_1.  \tag{3.2}
\end{align*}\\
By differentiating $\widehat{e^{sA}\varphi}$, we have 

\begin{align*}
\widehat{A\varphi}=-|\xi_1|\hat{\varphi} 
=i\frac{\xi_1}{|\xi_1|}(i\xi_1)\hat{\varphi}
=-\widehat{H\partial_1\varphi}.
 \tag{3.3}
\end{align*}\\
Since 

\begin{align*}
\widehat{e^{sA}H\varphi}=\widehat{Q_{s}*_1\varphi},  \tag{3.4}
\end{align*}\\
the Hilbert transform is represented by

\begin{align*}
H\varphi(x_1)=\lim_{s\to 0}(Q_s*_1\varphi)(x_1)
=\lim_{s\to 0}\frac{1}{\pi}\int_{\mathbb{R}}\frac{x_1-y_1}{(x_1-y_1)^{2}+s^{2}}\varphi(y_1)\dd y_1.    \tag{3.5}
\end{align*}\\
We use the kernels $P_s$ and $Q_s$ to calculate the Hilbert transform of $\nabla E$. Since 

\begin{equation*}
\begin{aligned}
\partial_{x_1}E(x)=-\frac{x_1}{2\pi(x_1^{2}+x_2^{2} )}
=-\frac{1}{2}Q_{x_2}(x_1),\\
\partial_{x_2}E(x)=-\frac{x_2}{2\pi(x_1^{2}+x_2^{2} )}
=-\frac{1}{2}P_{x_2}(x_1),\\
\end{aligned}
\tag{3.6}
\end{equation*}\\
and $\widehat{HP_s}=\widehat{Q_s}$, $s>0$, we have

\begin{equation*}
\begin{aligned}
H\partial_{x_2}E(x)
=
\begin{cases}
\quad \partial_{x_1} E(x),& \quad x_2>0,\\
\quad -\partial_{x_1} E(x),& \quad x_2<0.\\ 
\end{cases}
\end{aligned}
\tag{3.7}
\end{equation*}\\
To prove (1.22), we use the Hardy space ${\mathcal{H}}^{1}(\mathbb{R})$ \cite{Stein93}. See also \cite{Miyakawa96}. A tempered distribution $\varphi$ belongs to ${\mathcal{H}}^{1}(\mathbb{R})$ if 

\begin{align*}
\varphi^{+}(x_1)=\sup_{t>0}|(e^{t\partial_1^{2}}\varphi)(x_1)|\in L^{1}(\mathbb{R}).
\end{align*}\\
The quasi-norm for $\varphi\in {\mathcal{H}}^{1}(\mathbb{R})$ is defined by 

\begin{align*}
||\varphi||_{{\mathcal{H}}^{1}(\mathbb{R})}
=||\varphi^{+}||_{L^{1}(\mathbb{R})}.
\end{align*}\\
The space ${\mathcal{H}}^{1}(\mathbb{R})$ is smaller than $L^{1}(\mathbb{R})$. Indeed, $\varphi\in L^{1}(\mathbb{R})$ belongs to ${\mathcal{H}}^{1}(\mathbb{R})$ if and only if $H\varphi\in L^{1}(\mathbb{R})$ and the quasi-norm $||\cdot ||_{{\mathcal{H}}^{1}(\mathbb{R})}$ is equivalent to 

\begin{align*}
||\varphi||_{{\mathcal{H}}^{1}(\mathbb{R})}
\cong ||\varphi||_{L^{1}(\mathbb{R})}+||H\varphi ||_{L^{1}(\mathbb{R})}.
\end{align*}\\
If $H\varphi\in L^{1}(\mathbb{R})$, $\widehat{H\varphi}=-i\xi_1|\xi_1|^{-1}\hat{\varphi}(\xi_1)$ is continuous for $\xi_1\in \mathbb{R}$. Hence, $\varphi\in {\mathcal{H}}^{1}(\mathbb{R})$ implies 

\begin{align*}
0=\hat{\varphi}(0)=\int_{\mathbb{R}}\varphi(x_1)\dd x_1.
\end{align*}\\
Note that $\Gamma_0(x_1,1)=(4\pi)^{-1/2}e^{-|x_1|^{2}/4}\notin {\mathcal{H}}^{1}(\mathbb{R})$ since $\int_{\mathbb{R}}\Gamma_0(x_1,1)\dd x_1=1$. On the other hand, $\partial_1\Gamma_0(x_1,1)\in {\mathcal{H}}^{1}(\mathbb{R})$ since 

\begin{align*}
\partial_1\Gamma_0^{+}(x_1,1)
&=\sup_{t>0}|e^{t\partial_1^{2}}\partial_1\Gamma_0(x_1,1)|\\
&=\sup_{t>1}|\partial_1\Gamma_0(x_1,t)| \\
&=\sup\left\{\left|\frac{z_1}{x_1} \right|^{2} |\partial_1\Gamma_0(z_1,1)|\ \middle|\ |z_1|\leq |x_1|  \right\} 
\leq \frac{C}{1+|x_1|^{2}},\quad x_1\in \mathbb{R}.
\end{align*}\\
This implies 

\begin{align*}
A\Gamma_0(x_1,1)=-H\partial_1 \Gamma_0(x_1,1)\in L^{1}(\mathbb{R}).  \tag{3.8}
\end{align*}\\
We use (3.8) to prove (1.22) in Section 4.

\vspace{15pt}

\subsection{Solution formulas}
We prove the formulas (1.12) and (1.10). Let $u_0\in L^{2,\infty}_{\sigma}$ satisfy $\omega_0\in M$. Since $u_0=K\omega_0$ by Lemma 2.3 (i), it follows from (2.9) that

\begin{equation*}
\begin{aligned}
&S(t)u_0=\int_{\mathbb{R}^{2}_{+}}G(x,z,t)  \left(\int_{\mathbb{R}^{2}_{+}}\nabla_{z}^{\perp}D(z,y)\omega_0(\dd y) \right)\dd z
=\int_{\mathbb{R}^{2}_{+}} \left(\int_{\mathbb{R}^{2}_{+}}G(x,z,t)\nabla_{z}^{\perp}D(z,y)\dd z \right)\omega_0(\dd y).
\end{aligned}
\end{equation*}\\
By taking the rotation, we have

\begin{equation*}
\begin{aligned}
T(t)\omega_0&=-\nabla^{\perp} \cdot S(t)K\omega_0 
=\int_{\mathbb{R}^{2}_{+}}W(x,y,t)\omega_0(\dd y), \\  
W(x,y,t)&=-\nabla^{\perp}_{x}\cdot \left(\int_{\mathbb{R}^{2}_{+}}G(x,z,t)\nabla_{z}^{\perp}D(z,y)\dd z \right).
\end{aligned}
\tag{3.9}
\end{equation*}

\vspace{15pt}

%prop3.1
\begin{prop}
\begin{equation*}
\begin{aligned}
W(x,y,t)=&
-\partial_{x_2}\int_{\mathbb{R}^{2}_{+}}\left(\Gamma(x-z,t)-\Gamma(x-z^{*},t)\right)\partial_{z_2}D(z,y)\dd z \\
&-\partial_{x_1}\int_{\mathbb{R}^{2}_{+}}\left(\Gamma(x-z,t)-\Gamma(x-z^{*},t)\right)\partial_{z_1}D(z,y)\dd z  \\
&+2\int_{\mathbb{R}^{2}_{+} }\Gamma_0(x_2+z_2,t)(A\Gamma_0)(x_1-z_1,t)
\partial_{z_2}D(z,y)\dd z.
\end{aligned}
\tag{3.10}
\end{equation*}
\end{prop}

\vspace{15pt}

\begin{proof}
By (2.10),
\begin{equation*}
\begin{aligned}
G_{11}(x,z,t)&=\Gamma(x-z,t)-\Gamma(x-z^{*},t)-4\int_{0}^{x_2}\int_{\mathbb{R}}\partial_{x_1}^{2}E(x-w)\Gamma(w-z^{*},t)\dd w,\\
G_{12}(x,z,t)&= 0,\\
G_{21}(x,z,t)&=-4\int_{0}^{x_2}\int_{\mathbb{R}}\partial_{x_1}\partial_{x_2}E(x-w)\Gamma (w-z^{*},t)\dd w,\\
G_{22}(x,z,t)&=\Gamma(x-z,t)-\Gamma(x-z^{*},t).
\end{aligned}
\end{equation*}\\
It follows from (3.9) that

\begin{equation*}
\begin{aligned}
W(x,y,t)
=&\partial_{x_1}\left(\int_{\mathbb{R}^{2}_{+}} G_{21}(x,z,t)\partial_{z_2}D(z,y)\dd z-  \int_{\mathbb{R}^{2}_{+}} G_{22}(x,z,t)\partial_{z_1}D(z,y)\dd z \right)\\
 &-\partial_{x_2}\left(\int_{\mathbb{R}^{2}_{+}} G_{11}(x,z,t)\partial_{z_2}D(z,y)\dd z \right)\\
=&\int_{\mathbb{R}^{2}_{+}} \left(\partial_{x_1}G_{21}(x,z,t)-\partial_{x_2}G_{11}(x,z,t)   \right)\partial_{z_2}D(z,y)\dd z \\
&-  \int_{\mathbb{R}^{2}_{+}} \partial_{x_1}G_{22}(x,z,t)\partial_{z_1}D(z,y)\dd z. 
\end{aligned}
\tag{3.11}
\end{equation*}\\
Since 

\begin{align*}
\partial_{x_2}G_{11}(x,z,t)
&=\partial_{x_2}(\Gamma(x-z,t)-\Gamma(x-z^{*},t) ) 
 -4\int_{0}^{x_2}\int_{\mathbb{R}}\partial_{x_2}\partial_{x_1}^{2}E(x-w)\Gamma(w-z^{*},t)\dd w \\
&\quad -4\lim_{w_2\to x_2}\int_{\mathbb{R}}\partial_{x_1}^{2}E(x-w)\Gamma(w-z^{*},t)\dd w_1 \\
&=-\partial_{x_2}(\Gamma(x-z,t)-\Gamma(x-z^{*},t) ) 
+\partial_{x_1}G_{21}(x,z,t) \\
&\quad -4\lim_{w_2\to x_2}\int_{\mathbb{R}}\partial_{x_1}^{2}E(x-w)\Gamma(w-z^{*},t)\dd w_1,
\end{align*}\\
we have

\begin{align*}
\partial_{x_1}G_{21}(x,z,t)-\partial_{x_2}G_{11}(x,z,t)
=&-\partial_{x_2}(\Gamma(x-z,t)-\Gamma(x-z^{*},t) )  \\
&+4\lim_{w_2\to x_2}\int_{\mathbb{R}}\partial_{x_1}^{2}E(x-w)\Gamma(w-z^{*},t)\dd w_1.
\end{align*}\\
It follows from (3.11) that 

\begin{align*}
W(x,y,t)=
&-\int_{\mathbb{R}^{2}_{+}}\partial_{x_2}\left(\Gamma(x-z,t)-\Gamma(x-z^{*},t)\right)\partial_{z_2}D(z,y)\dd z \\
&-\int_{\mathbb{R}^{2}_{+}}\partial_{x_1}(\Gamma(x-z,t)-\Gamma(x-z^{*},t))\partial_{z_2}D(z,y)\dd z\\
&+4\int_{\mathbb{R}^{2}_{+}}\left(\lim_{w_2\to x_2}\int_{\mathbb{R}}\partial_{x_1}^{2}E(x-w)\Gamma(w-z^{*},t)\dd w_1\right)\partial_{z_1}D(z,y)\dd z
\end{align*}\\
To prove (3.10), it suffices to show that 

\begin{align*}
\lim_{w_2\to x_2}\int_{\mathbb{R}}\partial_{x_1}^{2}E(x-w)\Gamma(w-z^{*},t)\dd w_1
=\frac{1}{2}\Gamma_0(x_2+z_2,t)(A\Gamma_0)(x_1-z_1,t).  \tag{3.12}
\end{align*}\\
Since $\partial_{x_1}E(x)=-2^{-1}Q_{x_2}(x_1)$ by (3.6) and $\Gamma(x,t)=\Gamma_0(x_1,t)\Gamma_0(x_2,t)$, integration by parts yields 

\begin{align*}
\int_{\mathbb{R}}\partial_{x_1}^{2}E(x-w)\Gamma(w-z^{*},t)\dd w_1
&=-\int_{\mathbb{R}}\partial_{w_1}\partial_{x_1}E(x-w)\Gamma(w-z^{*},t)\dd w_1 \\
&=\int_{\mathbb{R}}\partial_{x_1}E(x-w)\partial_{w_1}\Gamma(w-z^{*},t)\dd w_1\\
&=-\frac{1}{2}\Gamma_0(w_2+z_2,t)\int_{\mathbb{R}}Q_{x_2-w_2}(x_1-w_1)\partial_{w_1}\Gamma_{0}(w_1-z_1,t)\dd w_1 \\
&=-\frac{1}{2}\Gamma_0(w_2+z_2,t)(Q_{x_2-w_2}*_1\partial_{1}\Gamma_{0})(x_1-z_1,t).
\end{align*}\\
By (3.5) and (3.3), we have

\begin{align*}
(Q_{x_2-w_2}*_1\partial_{1}\Gamma_{0})(x_1-z_1,t)
\to (H\partial_1\Gamma_0)(x_1-z_1,t)=-(A\Gamma_0)(x_1-z_1,t)\quad \textrm{as}\ w_2\to x_2.
\end{align*}\\
We proved (3.12). The proof is complete. 
\end{proof}

\vspace{15pt}

%prop2.2
\begin{prop}
Set $e^{t\Delta_{N}}$, $e^{t\Delta_{D}}$ and $e^{t\partial^{2}_{1}}$ by (1.11). Then, 

\begin{align*}
e^{t\Delta_{D}}\partial_2 \varphi&=\partial_2e^{t\Delta_{N}}\varphi,  \tag{3.13}\\
e^{t\Delta_{N}}\partial_2 \varphi&=\partial_2e^{t\Delta_{D}}\varphi-2\Gamma_0(x_2,t)e^{t\partial_1^{2}}\varphi(\cdot,0). \tag{3.14} 
\end{align*}
\end{prop}

\vspace{5pt}

\begin{proof}
By integration by parts, 

\begin{align*}
e^{t\Delta_N}\partial_2 \varphi
&=\int_{\mathbb{R}^{2}_{+}}(\Gamma(x-y,t)+\Gamma(x-y^{*},t) )\partial_{y_2}\varphi(y)\dd y \\
&=-\int_{\mathbb{R}^{2}_{+}}\partial_{y_2}(\Gamma(x-y,t)+\Gamma(x-y^{*},t) )\varphi(y)\dd y+\left[\int_{\mathbb{R}}(\Gamma(x-y,t)+\Gamma(x-y^{*},t))\varphi(y)\dd y_1   \right]_{y_2=0}^{y_2=\infty}\\
&=\partial_{x_2}\int_{\mathbb{R}^{2}_{+}}(\Gamma(x-y,t)-\Gamma(x-y^{*},t) )\varphi(y)\dd y
-2\Gamma_0(x_2,t)\int_{\mathbb{R}}\Gamma_0(x_1-y_1,t)\varphi(y_1,0)\dd y_1.
\end{align*}\\
Thus (3.14) holds. In a similar way, (3.13) follows.
\end{proof}

\vspace{15pt}

%lem 2.3
\begin{lem}
The formulas (1.10) and (1.12) hold for $u_0\in L^{2,\infty}_{\sigma}$ satisfying $\omega_0\in M$.
\end{lem}

\vspace{5pt}

\begin{proof}
We substitute $\varphi(z_1,z_2)=\partial_{z_2}D(z,y)$ into

\begin{align*}
e^{t\Delta_N}\partial_2\varphi
=\int_{\mathbb{R}^{2}_{+}}\left(\Gamma(x-z,t)+\Gamma(x-z^{*},t) \right)\partial_{z_2}\varphi(z)\dd z.
\end{align*}\\
Since

\begin{align*}
\partial_{z_2}D(z,y)=\partial_{z_2}\left(E(z-y)-E(z-y^{*})\right)
=-\frac{1}{2\pi}\left( \frac{z_2-y_2}{|z-y |^{2}}-\frac{z_2+y_2}{|z-y^{*} |^{2}} \right),
\end{align*}\\
we observe that 

\begin{align*}
\varphi(z_1,0)=\frac{y_2}{\pi (|z_1-y_1|^{2}+y_2^{2} )}=P_{y_2}(z_1-y_1).
\end{align*}\\
It follows from (3.14) that

\begin{align*}
e^{t\Delta_N}\partial_{z_2}^{2}D=\partial_{x_2}e^{t\Delta_D}\partial_{z_2}D
-2\Gamma_0(x_2,t)\int_{\mathbb{R}}\Gamma_0(x_1-z_1,t)P_{y_2}(z_1-y_1)\dd z_1.  \tag{3.15}
\end{align*}\\
By (3.10) and (3.15), we have

\begin{align*}
W(x,y,t)&=-\partial_{x_2}e^{t\Delta_{D}}\partial_{z_2}D-\partial_{x_1}e^{t\Delta_{D}}\partial_{z_1}D+2\int_{\mathbb{R}^{2}_{+} }\Gamma_0(x_2+z_2,t)(A\Gamma_0)(x_1-z_1,t)
\partial_{z_2}D(z,y)\dd z\\
&=-e^{t\Delta_N}\partial_{z_2}^{2}D-2\Gamma_0(x_2,t)\int_{\mathbb{R}}\Gamma_0(x_1-z_1,t)P_{y_2}(z_1-y_1)\dd z_1-e^{t\Delta_{D}}\partial_{z_1}^{2}D \\
&\quad +2\int_{\mathbb{R}^{2}_{+} }\Gamma_0(x_2+z_2,t)(A\Gamma_0)(x_1-z_1,t)\partial_{z_2}D(z,y)\dd z.
\end{align*}\\
Since $-\Delta_{z} D(z,y)=\delta_{y}$ with the Dirac measure $\delta_{y}$ in $\overline{\mathbb{R}^{2}_{+}}$ at $z=y$, we have

\begin{equation*}
\begin{aligned}
W(x,y,t)&=e^{t\Delta_{N}}\delta_{y}+(e^{t\Delta_N}-e^{t\Delta_D})\partial_{z_1}^{2}D-2\Gamma_0(x_2,t)\int_{\mathbb{R}}\Gamma_0(x_1-z_1,t)P_{y_2}(z_1-y_1)\dd z_1\\
&\quad +2\int_{\mathbb{R}^{2}_{+} }\Gamma_0(x_2+z_2,t)(A\Gamma_0)(x_1-z_1,t)\partial_{z_2}D(z,y)\dd z.
\end{aligned}
\tag{3.16}
\end{equation*}\\
Since  

\begin{align*}
{\mathcal{F}}(A\Gamma_0*_1\partial_{z_2}D)
=-|\xi_1|\widehat{\Gamma_0}\widehat{\partial_{z_2}D}
={\mathcal{F}}(\Gamma_0*_1A\partial_{z_2}D),
\end{align*}\\
we have 

\begin{align*}
\int_{\mathbb{R} }(A\Gamma_0)(x_1-z_1,t)\partial_{z_2}D(z,y)\dd z_1
=\int_{\mathbb{R}}\Gamma_0(x_1-z_1,t)(A\partial_{z_2}D)(z,y)\dd z_1.
\end{align*}\\
Hence, 

\begin{align*}
2\int_{\mathbb{R}^{2}_{+} }\Gamma_0(x_2+z_2,t)(A\Gamma_0)(x_1-z_1,t)\partial_{z_2}D(z,y)\dd z
&=2\int_{\mathbb{R}^{2}_{+} }\Gamma_0(x_2+z_2,t)\Gamma_0(x_1-z_1,t)(A\partial_{z_2}D)(z,y)\dd z \\
&=2\int_{\mathbb{R}^{2}_{+} }\Gamma(x-z^{*},t)(A\partial_{z_2}D)(z,y)\dd z\\
&=(e^{t\Delta_N}-e^{t\Delta_D})A\partial_{z_2}D.
\end{align*}\\
Since $A=-H\partial_1$ by (3.3), this implies 

\begin{align*}
&(e^{t\Delta_N}-e^{t\Delta_D})\partial_{z_1}^{2}D
+2\int_{\mathbb{R}^{2}_{+} }\Gamma_0(x_2+z_2,t)(A\Gamma_0)(x_1-z_1,t)\partial_{z_2}D(z,y)\dd z\\
&=(e^{t\Delta_N}-e^{t\Delta_D})(\partial_{z_1}^{2}D+A\partial_{z_2}D) \\
&=-(e^{t\Delta_N}-e^{t\Delta_D})(H\partial_{z_2}-\partial_{z_1})\partial_{z_1}D.
\end{align*}\\
It follows from (3.16) that 

\begin{equation*}
\begin{aligned}
W(x,y,t)=
&e^{t\Delta_N}\delta_y-(e^{t\Delta_N}-e^{t\Delta_D})(H\partial_{z_2}-\partial_{z_1})\partial_{z_1}D \\
&-2\Gamma_0(x_2,t)\int_{\mathbb{R}}\Gamma_0(x_1-z_1,t)P_{y_2}(z_1-y_1)\dd z_1.     
\end{aligned}
\tag{3.17}
\end{equation*}\\
Since 

\begin{align*}
u^{1}_{0}(z_1,0)
=\frac{1}{\pi}\int_{\mathbb{R}^{2}_{+}}\frac{y_2}{|z_1-y_1|^{2}+y_2^{2} }\omega_0(\dd y)
=\int_{\mathbb{R}^{2}_{+}}P_{y_2}(z_1-y_1)\omega_0(\dd y),
\end{align*}\\
by (1.2), we have 

\begin{align*}
\int_{\mathbb{R}^{2}_{+}}\left(\int_{\mathbb{R}}\Gamma_0(x_1-z_1,t)P_{y_2}(z_1-y_1)\dd z_1 \right)\omega_0(\dd y)
&=\int_{\mathbb{R}}\Gamma_0(x_1-z_1,t) \left(\int_{\mathbb{R}^{2}_{+}}P_{y_2}(z_1-y_1)\omega_0(\dd y)\right)\dd z_1 \\
&=\int_{\mathbb{R}}\Gamma_0(x_1-z_1,t)u^{1}_{0}(z_1,0)\dd z_1
=e^{t\partial_1^{2}}u^{1}(\cdot,0).
\end{align*}\\
Thus, (1.12) follows by integrating (3.17) by the measure $\omega_0\in M$.

To prove (1.10), we set 

\begin{align*}
R(z)&=(H\partial_{z_2}-\partial_{z_1})\partial_{z_1}E(z),\\
\Pi(z,y)&=(H\partial_{z_2}-\partial_{z_1})\partial_{z_1}D(z,y)=R(z-y)-R(z-y^{*}).
\end{align*}\\
By (3.7),

\begin{align*}
R(z)=\partial_{z_1} \left(H\partial_{z_2} E(z)-\partial_{z_1} E(z)\right)
=
\begin{cases}
\quad 0, & z_2>0, \\
\quad -2\partial_{z_1}^{2}E(z), & z_2<0.
\end{cases}
\end{align*}\\
This implies $R(z-y^{*})=0$ for $z, y\in \mathbb{R}^{2}_{+}$ and 

\begin{align*}
\Pi(z,y)=R(z-y)=
\begin{cases}
\quad 0,& y_2<z_2,\\
\quad -2\partial_{z_1}^{2}E(z-y), & 0<z_2<y_2.
\end{cases} 
\end{align*}\\
Hence

\begin{align*}
(e^{t\Delta_N}-e^{t\Delta_D} )(H\partial_{z_2}-\partial_{z_1})\partial_{z_1}D
&=2\int_{\mathbb{R}^{2}_{+}}\Gamma(x-z^{*},t)\Pi(z,y)\dd z \\
&=-4\int_{0}^{y_2}\int_{\mathbb{R}}\Gamma(x-z^{*},t)\partial_{z_1}^{2}E(z-y)\dd z. 
\end{align*}\\
By (3.17), (1.10) follows.
\end{proof}

\vspace{5pt}

\subsection{A kernel estimate}

We give a pointwise estimate for $W(x,y,t)$.

\vspace{15pt}

%lem2.4
\begin{lem}
(i) Set

\begin{equation*}
\begin{aligned}
W(x,y,t)&=\Gamma(x-y,t)+W^{*}(x,y,t)+W_{\textrm{tr}}(x,y,t),\\
W^{*}(x,y,t)&=\Gamma(x-y^{*},t)+4\int_{0}^{y_2}\int_{\mathbb{R}}\Gamma (x-z^{*}, t)\partial_{z_1}^{2}E(z-y)\dd z,\\
W_{\textrm{tr}}(x,y,t)&=-2\Gamma_0(x_2,t)\int_{\mathbb{R}}\Gamma_0(x_1-z_1,t)P_{y_2}(z_1-y_1)\dd z_1.
\end{aligned}
\tag{3.18}
\end{equation*}\\
Then, 

\begin{align*}
W(x,y,t)&=\lambda^{2}W(\lambda x,\lambda y,\lambda^{2} t),\quad \lambda>0,  \tag{3.19} \\
W^{*}(x,y,t)&=-G^{*}_{11}(y,x,t),  \tag{3.20} \\
|\partial_t^{s}\partial_x^{k}\partial_y^{m}W^{*}(x,y,t)|
&\leq \frac{Ce^{-cx_2^{2}/t}}{t^{s+k_2/2}(y_2^{2}+t)^{m_2/2}(|x-y^{*}|^{2}+t)^{(k_1+m_1+2)/2}  }, \tag{3.21}
\end{align*}\\
for $k=(k_1,k_2)$, $m=(m_1,m_2)$, and $s\geq 0$.

\noindent 
(ii) Set 

\begin{align*}
\tilde{W}(x,y,t)=4\int_{0}^{y_2}\int_{\mathbb{R}}\Gamma (x-z^{*}, t)\partial_{z_1}^{2}E(z-y)\dd z.  \tag{3.22}
\end{align*}\\
Then, 

\begin{align*}
-(e^{t\Delta_N}-e^{t\Delta_D})(H\partial_2-\partial_1 )\partial_1(-\Delta_D)^{-1}\omega_0=\int_{\mathbb{R}^{2}_{+}}\tilde{W}(x,y,t)\omega_0(\dd y),\quad \omega_0\in M.  \tag{3.23}
\end{align*}
\end{lem}

\vspace{5pt}

\begin{proof}
The functions $\Gamma(x,t)$, $\nabla^{2}E(x)$ and $P_{s}(x_1)$ satisfy the scaling properties

\begin{align*}
\Gamma(x,t)&=\lambda^{2} \Gamma( \lambda x, \lambda^{2} t),\\
\nabla^{2} E(x)&=\lambda^{2} (\nabla^{2}E)( \lambda x),\\
P_{s}(x_1)&=\lambda P_{\lambda s}(\lambda x_1),\quad \lambda>0.
\end{align*}\\
By the changing variable, we observe that 

\begin{align*}
W^{*}(\lambda x,\lambda y,\lambda^{2} t)
&=\Gamma(\lambda(x-y^{*}),\lambda^{2}t  )
+4\int_{0}^{\lambda y_2}\int_{\mathbb{R}}\Gamma(\lambda x-z^{*},\lambda^{2}t)\partial_{z_1}^{2}E(z-\lambda y)\dd z \\
&=\Gamma(\lambda(x-y^{*}),\lambda^{2}t  )
+4\int_{0}^{ y_2}\int_{\mathbb{R}}\Gamma(\lambda (x-w^{*}),\lambda^{2}t)(\partial_{1}^{2}E)(\lambda(w- y))\lambda^{2}\dd w\\
&=\lambda^{-2}\left(\Gamma(x-y^{*},t)+4 \int_{0}^{ y_2}\int_{\mathbb{R}}\Gamma(x-w^{*},t)\partial_{w_1}^{2}E(w- y)\dd w\right)\\
&=\lambda^{-2}W(x,y,t),
\end{align*}
\begin{align*}
W_{\textrm{tr}}(\lambda x,\lambda y,\lambda^{2} t)
&=-2\Gamma_{0}(\lambda x_2,\lambda^{2}t )\int_{\mathbb{R}}\Gamma_0(\lambda x_1-z_1,\lambda^{2} t)P_{\lambda y_2}(z_1-\lambda y_1)\dd z_1\\
&=-2 \Gamma_0(\lambda x_2,\lambda^{2} t)
\int_{\mathbb{R}}\Gamma_0(\lambda (x_1-w_1),\lambda^{2} t)P_{\lambda y_2}(\lambda(w_1- y_1))\lambda \dd w_1\\
&=\lambda^{-2}W_{\textrm{tr}}(x, y, t).
\end{align*}\\
Thus (3.19) holds. To prove (3.20), we observe from (2.10) that 

\begin{align*}
G^{*}_{11}(x,z,t)=-\Gamma(x-z^{*},t)-4\int_{0}^{x_2}\int_{\mathbb{R}}\partial_{x_1}^{2}E(x-w)\Gamma(w-z^{*},t)\dd w.
\end{align*}\\
By replacing $x$ and $z=y$, we have 

\begin{align*}
G^{*}_{11}(y,x,t)=-\Gamma(y-x^{*},t)-4\int_{0}^{y_2}\int_{\mathbb{R}}\partial_{y_1}^{2}E(y-w)\Gamma(w-x^{*},t)\dd w.
\end{align*}\\
We change the variable $w$ to $z$. Since $\Gamma(y-x^{*},t)=\Gamma(x-y^{*},t)$ and 

\begin{align*}
\partial_{y_1}^{2}E(y-z)=
-\frac{|y_2-z_2|^{2}-|y_1-z_1|^{2}}{2\pi |y-z|^{4}}
=\partial_{z_1}^{2}E(z-y),
\end{align*}\\
it follows from (3.18) that

\begin{align*}
G^{*}_{11}(y,x,t)
&=-\Gamma(y-x^{*},t)-4\int_{0}^{y_2}\int_{\mathbb{R}}\partial_{y_1}^{2}E(y-z)\Gamma(z-x^{*},t)\dd z \\
&=-\Gamma(x-y^{*},t)-4\int_{0}^{y_2}\int_{\mathbb{R}}\partial_{z_1}^{2}E(z-y)\Gamma(x-z^{*},t)\dd z\\
&=-W^{*}(x,y,t).
\end{align*}\\
We proved (3.20). The estimate (3.21) follows from (3.20) and (2.11). Since $W(x,y,t)=\Gamma(x-y,t)+\Gamma(x-y^{*},t)+\tilde{W}(x,y,t)+W_{\textrm{tr}}(x,y,t)$, (3.23) follows from (1.12).
\end{proof}

\vspace{15pt}

%rem 2.4
\begin{rems}
(i) The operator $T(t): \omega_0\longmapsto \omega(\cdot,t)$ is a solution operator of the heat equation,

\begin{equation*}
\begin{aligned}
\partial_t\omega-\Delta \omega&=0\qquad \textrm{in}\ \Omega\times (0,\infty),\\
\frac{\partial \omega}{\partial n}+A\omega&=0\qquad \textrm{on}\ \partial\Omega\times (0,\infty),\\
\omega&=\omega_0\hspace{15pt} \textrm{on}\ \Omega\times \{t=0\},
\end{aligned}
\tag{3.24}
\end{equation*}\\
for the half plane $\Omega=\mathbb{R}^{2}_{+}$, where $\partial /\partial n=n\cdot \nabla$ denotes the normal derivative and $n={}^{t}(0,-1)$. The equation (3.24) is the vorticity equations associated with the Stokes equations

\begin{equation*}
\begin{aligned}
\partial_t v-\Delta v+\nabla q=0,\ \D\ v&=0\qquad  \textrm{in}\ \Omega\times (0,\infty),\\
v &=0\qquad \textrm{on}\ \partial\Omega\times (0,\infty),\\
v&=v_0\hspace{18pt} \textrm{on}\ \Omega\times \{t=0\}.
\end{aligned}
\tag{3.25}
\end{equation*}\\
Since $-\Delta v=\nabla^{\perp}\omega$ and the pressure $q$ satisfies the Neumann problem

\begin{align*}
-\Delta q&=0\quad \textrm{in}\ \mathbb{R}^{2}_{+},\\
\frac{\partial q}{\partial x_2}&=\Delta v^{2}=\partial_1 \omega\quad \textrm{on}\ \partial \mathbb{R}^{2}_{+},
\end{align*}\\
$q$ is represented as

\begin{align*}
q=-\int_{x_2}^{\infty}e^{sA}\partial_1 \omega\dd s
\end{align*}\\
by the Poisson semigroup $e^{sA}$ defined by (3.2). Hence,

\begin{align*}
\partial_1 q
=-\int_{x_2}^{\infty}e^{sA}\partial_1^{2}\omega\dd s
=\int_{x_2}^{\infty}e^{sA}A^{2}\omega\dd s
=-e^{x_2A}A \omega.
\end{align*}\\
Thus taking the tangential trace to $(3.25)_1$ implies 

\begin{align*}
0=\lim_{x_2\to0} (\partial_tv^{1}+\partial_2\omega+\partial_1 q)=(\partial_2-A)\omega.
\end{align*}\\

\noindent 
(ii) The formula (1.12) gives a solution to (3.24). Since $A=-H\partial_1$ by (3.3), (1.12) is represented by

\begin{align*}
T(t)\omega_0=e^{t\Delta_N}\omega_0
+(e^{t\Delta_N}-e^{t\Delta_D})A(\partial_2-A)(-\Delta_D)^{-1}\omega_0
-2\Gamma_0(x_2,t)e^{t\partial_1^{2}}u^{1}_{0}(\cdot,0).  \tag{3.26}
\end{align*}\\
By (3.13) and (3.14), we have

\begin{align*}
(\partial_2-A)(e^{t\Delta_N}-e^{t\Delta_D})\varphi
=-(e^{t\Delta_N}-e^{t\Delta_D})(\partial_2+A)\varphi-2\Gamma_0(x_2,t)e^{t\partial_1^{2}}\varphi(\cdot,0).
\end{align*}\\
We substitute $\varphi=A(\partial_2-A)(-\Delta_D)^{-1}\omega_0$. Since $(\partial_2+A)\varphi=-A \omega_0$ and 

\begin{align*}
\varphi(x_1,0)=A\partial_2(-\Delta_D)^{-1}\omega_0+\partial_1^{2}(-\Delta_D)^{-1}\omega_0 
=Au^{1}_{0}(x_1,0)
+\partial_1^{2}u^{2}_{0}(x_1,0)=Au^{1}_{0}(x_1,0)
\end{align*}\\
by (1.2), it follows that

\begin{align*}
(\partial_2-A)(e^{t\Delta_N}-e^{t\Delta_D})A(\partial_2-A)(-\Delta_D)^{-1}\omega_0
=(e^{t\Delta_N}-e^{t\Delta_D})A \omega_0
-2\Gamma_0(x_2,t)e^{t\partial_1^{2}}A u^{1}_{0}(\cdot,0).
\end{align*}\\
Multiplying $(\partial_2-A)$ by (3.26) implies 

\begin{align*}
(\partial_2-A)T(t)\omega_0=
&(\partial_2-A)e^{t\Delta_N}\omega_0
+(e^{t\Delta_N}-e^{t\Delta_D})A\omega_0
-2\Gamma_0(x_2,t)e^{t\partial_1^{2}}A u^{1}_{0}(\cdot,0) \\
&-2(\partial_{x_2}-A)(\Gamma_0(x_2,t)e^{t\partial_1^{2}} u^{1}_{0}(\cdot,0))\\
=&(\partial_2-A)e^{t\Delta_N}\omega_0
+A(e^{t\Delta_N}-e^{t\Delta_D})\omega_0
-2\partial_{x_2}\Gamma_0(x_2,t)e^{t\partial_1^{2}} u^{1}_{0}(\cdot,0). 
\end{align*}\\
Since $\partial_{x_2}\Gamma_0(0,t)=0$, sending $x_2\to 0$ yields $(\partial_2-A)T(t)\omega_0=0$ for $x_2=0$. We prove the convergence to initial data (1.15) in Lemma 4.1.

\noindent
(iii) We are able to write the vorticity equations by (3.24) even for domains $\Omega$ by using the operators $A=-H\partial_{\textrm{tan}}$ and $H: g\longmapsto -\partial_{\textrm{tan}} q$, associated with the Neumann problem

\begin{align*}
-\Delta q&=0\quad \textrm{in}\ \Omega,\\
\frac{\partial q}{\partial n} &=g\quad \textrm{on}\ \partial\Omega.
\end{align*}\\
Here, $\partial_{\textrm{tan}}=n\cdot\nabla^{\perp}=n^{\perp}\cdot \nabla$ for the unit outward normal vector field $n={}^{t}(n^{1},n^{2})$ on $\partial\Omega$ and $n^{\perp}={}^{t}(-n^{2},n^{1})$. Since $\D\ v=0$ and $-\Delta v=\nabla^{\perp}\omega$, the pressure $q$ solves the Neumann problem for $g=\Delta v\cdot n=- \partial_{\textrm{tan}}\omega$ \cite{JL04}, \cite{LLP}. Hence,

\begin{align*}
-\partial_{\textrm{tan}}q=Hg=-H\partial_{\textrm{tan}}\omega=A\omega.
\end{align*}\\
Thus multiplying $-n^{\perp}$ by $(3.25)_1$ and taking the trace implies

\begin{align*}
0=-n^{\perp}\cdot \nabla^{\perp}\omega-n^{\perp}\cdot\nabla q 
= \frac{\partial \omega}{\partial n}-\partial_{\textrm{tan}}q
=\frac{\partial \omega}{\partial n}+A\omega.
\end{align*}\\
The operator $T(t)=-\nabla^{\perp}\cdot S(t)K$ can be defined also for domains. For example, if $\Omega$ is bounded and simply-connected, the Biot-Savart law $u_0=K\omega_0$ is available and we are able to define $T(t)$ in the same way as the half plane. It is an interesting question whether the $L^{1}$-estimate (1.13) holds for domains. The $L^{\infty}$-estimate (1.14) is still valid for bounded domains \cite{AG1}, while the $L^{1}$-boundedness of $S(t)$ has been an open question \cite[Remark 5.2]{DHP}.
\end{rems}

\vspace{15pt}

%section3
\section{The semigroup associated with vorticity}

\vspace{15pt}

We study continuity of $T(t)\omega_0$ as $t\to 0$ by using the formula (1.12) for (i) general $\omega_0\in M$, (ii) continuous measures $\omega_0\in M$ and (iii) $\omega_0\in L^{1}$ satisfying $u^{1}_{0}(x_1,0)=0$. In the last subsection, we prove the asymptotic formula (1.22) as $t\to\infty$.

\vspace{15pt}

\subsection{Continuity in the vague topology}

We shall show that $T(t)$ forms a (not strongly continuous) bounded analytic semigroup on $M$.

\vspace{15pt}

%lem3.1
\begin{lem}
\begin{align*}
||\partial_{t}^{s}\partial_{x}^{k}T(t)\omega_0||_{q}\leq \frac{C}{t^{1-1/q+|k|/2+s}}||\omega_0||_{M},\quad t>0,\ \omega_0\in M,\ 1\leq q\leq \infty.   \tag{4.1} 
\end{align*}
\begin{align*}
T(t+s)=T(t)T(s)\quad t,s\geq 0.    \tag{4.2}  
\end{align*}
\begin{align*}
T(t)\omega_0\to \omega_0-\delta_{0,[0,\infty)}u^{1}_{0}(x_1,0)\quad \textrm{varguely on}\ M\quad \textrm{as}\ t\to0.   \tag{4.3}
\end{align*}\\
Set $u_0=K\omega_0$ and $\tilde{\omega}_0=\omega_0-\delta_{0,[0,\infty)}u^{1}_{0}(\cdot,0)$ for $\omega_0\in M$. Then,

\begin{align*}
T(t)\tilde{\omega}_0=T(t)\omega_0\to \tilde{\omega}_0\quad \textrm{varguely on}\ M\quad \textrm{as}\ t\to0.   \tag{4.4}
\end{align*}

\end{lem}

\vspace{5pt}

\begin{proof}
We show (4.1) for $q=1$, $s=0$ and $k=0$. The case $q\in (1,\infty]$, $s\geq 0$ and $|k|\geq 0$ follows in the same way. By (1.12),

\begin{align*}
T(t)\omega_0
=e^{t\Delta_N}\omega_0-(e^{t\Delta_N}-e^{t\Delta_D})(H\partial_2-\partial_1)\partial_1(-\Delta_{D})^{-1}\omega_0-2\Gamma_0(x_2,t)e^{t\partial_1^{2}}u^{1}_{0}(\cdot,0).
\end{align*}\\
By $(2.7)_1$, we have $||u^{1}_{0}||_{L^{1}(\mathbb{R})}\leq ||\omega_0||_{M}$. By estimating the kernels for $e^{t\Delta_N}$ and $e^{t\partial_1^{2}}$ defined by (1.11), we see that the first and third terms satisfy the desired estimate. We estimate the second term. By (3.23) and (3.21), we observe that 

\begin{align*}
|(e^{t\Delta_N}-e^{t\Delta_D})(H\partial_2-\partial_1)\partial_1(-\Delta_D)^{-1}\omega_0|
\lesssim \int_{\mathbb{R}^{2}_{+}}\frac{e^{-cx_2^{2}/t}}{|x-y^{*}|^{2}+t}|\omega_0|(\dd y).
\end{align*}\\
Since 

\begin{align*}
\int_{\mathbb{R}^{2}_{+}}\frac{e^{-cx_2^{2}/t}}{|x-y^{*}|^{2}+t}\dd x
=\int_{0}^{\infty}e^{-cx_2^{2}/t}\left(\int_{\mathbb{R}}\frac{\dd z_1}{z_1^{2}+|x_2+y_2|^{2}+t}\right)\dd x_2 
&=C\int_{0}^{\infty}\frac{e^{-cx_2^{2}/t} }{(|x_2+y_2|^{2}+t)^{1/2}}\dd x_2 \\
&\leq C\int_{0}^{\infty}\frac{e^{-cx_2^{2}/t} }{(|x_2|^{2}+t)^{1/2}}\dd x_2\\
&=C',
\end{align*}\\
we obtain 

\begin{align*}
||(e^{t\Delta_N}-e^{t\Delta_D})(H\partial_2-\partial_1)\partial_1(-\Delta_D)^{-1}\omega_0||_{1}\leq C||\omega_0||_{M},\quad t>0.
\end{align*}\\
We proved (4.1). Since $u=K\omega=K(-\nabla^{\perp})\cdot u$ for $u\in L^{2,\infty}_{\sigma}$ with $\omega_0\in M$ by Lemma 2.3 (i), 

\begin{align*}
T(t+s)=-\nabla^{\perp}\cdot S(t+s)K
=-\nabla^{\perp}\cdot S(t)K(-\nabla^{\perp})\cdot S(s)K
=T(t)T(s).
\end{align*}\\
Thus (4.2) holds. To prove (4.3), we take $\varphi\in C_{0}$ and set 

\begin{align*}
(T(t)\omega_0,\varphi)
&=(e^{t\Delta_N}\omega_0,\varphi)
-\left((e^{t\Delta_N}-e^{t\Delta_D})(H\partial_2-\partial_1)\partial_1(-\Delta_{D})^{-1}\omega_0,\varphi \right)
-\left(2\Gamma_0(\cdot,t)e^{t\partial_1^{2}}u^{1}_{0}(\cdot,0),\varphi \right) \\
&=I(t)+II(t)+III(t),
\end{align*}\\
with the pairing $(\cdot,\cdot)$ for $M$ and $C_0$. Let $\varphi^{\textrm{even}}$ be the even extension of $\varphi$, i.e., 

\begin{equation*}
\begin{aligned}
\varphi^{\textrm{even}}(x)=
\begin{cases}
&\varphi(x_1,x_2),\quad x_2\geq 0,\\
&\varphi(x_1,-x_2),\quad x_2< 0.
\end{cases}
\end{aligned}
\tag{4.5}
\end{equation*}\\
Observe that 

\begin{align*}
e^{t\Delta_{N}}\varphi
=\int_{\mathbb{R}^{2}_{+}}\left(\Gamma(x-y,t)+\Gamma(x-y^{*},t)\right)\varphi(y)\dd y
=\int_{\mathbb{R}^{2}}\Gamma(x-y,t)\varphi^{\textrm{even}}(y)\dd y
=:e^{t\Delta}\varphi^{\textrm{even}}.
\end{align*}\\
Since $\varphi^{\textrm{even}}\in C_{0}(\mathbb{R}^{2})$, it follows that 

\begin{align*}
||e^{t\Delta_{N}}\varphi-\varphi||_{L^{\infty}(\mathbb{R}^{2}_{+})}
\leq ||e^{t\Delta }\varphi^{\textrm{even}}-\varphi^{\textrm{even}}||_{L^{\infty}(\mathbb{R}^{2})}\to 0\quad \textrm{as}\ t\to 0.
\end{align*}\\
By the Fubini's theorem, 

\begin{align*}
I(t)=(e^{t\Delta_N}\omega_0,\varphi)
&=\int_{\mathbb{R}^{2}_{+}}\left( \int_{\mathbb{R}^{2}_{+}}\left(\Gamma(x-y,t)+\Gamma(x-y^{*},t) \right)\omega_0(\dd y) \right)\varphi(x)\dd x \\
&=\int_{\mathbb{R}^{2}_{+}}\left( \int_{\mathbb{R}^{2}_{+}}\left(\Gamma(x-y,t)+\Gamma(x-y^{*},t) \right)\varphi(x)\dd x \right)\omega_0(\dd y) \\
&=\int_{\mathbb{R}^{2}_{+}}\left( \int_{\mathbb{R}^{2}_{+}}\left(\Gamma(y-x,t)+\Gamma(y-x^{*},t) \right)\varphi(x)\dd x \right)\omega_0(\dd y) \\
&=(\omega_0, e^{t\Delta_N}\varphi )\to (\omega_0,\varphi)\quad \textrm{as}\ t\to 0.
\end{align*}\\
By (3.23),

\begin{align*}
II(t)=
\int_{\mathbb{R}^{2}_{+}}\left(\int_{\mathbb{R}^{2}_{+}}\tilde{W}(x,y,t)\omega_0(\dd y)   \right)\varphi(x)\dd x
=
\int_{\mathbb{R}^{2}_{+}}\left(\int_{\mathbb{R}^{2}_{+}}\tilde{W}(x,y,t)\varphi(x)\dd x \right)\omega_0(\dd y).
\end{align*}\\ 
Since $\tilde{W}(x,y,t)$ satisfies (3.21), we estimate

\begin{align*}
\left|\int_{\mathbb{R}^{2}_{+}}\tilde{W}(x,y,t)\varphi(x)\dd x \right|
\lesssim \int_{\mathbb{R}^{2}_{+}}\frac{e^{-cx_2^{2}/t }}{(|x-y^{*}|^{2}+t)}|\varphi(x)|\dd x
&\lesssim \int_{0}^{\infty}\frac{e^{-cx_2^{2}/t}}{((x_2+y_2)^{2}+t)^{1/2}}
||\varphi||_{L^{\infty}(\mathbb{R})}(x_2)\dd x_2 \\
&=:\rho(y_2,t).
\end{align*}\\
Observe that $\sup_{t>0}||\rho||_{\infty}(t)\leq C$ for some $C>0$ and $\lim_{t\to 0}\rho(y_2,t)=0$ for each $y_2>0$. The dominated convergence theorem yields

\begin{align*}
|II(t)|
=\left|\int_{\mathbb{R}^{2}_{+}}\left(\int_{\mathbb{R}^{2}_{+}}\tilde{W}(x,y,t)\varphi(x)\dd x \right)\omega_0(\dd y)\right| 
\lesssim \int_{\mathbb{R}^{2}_{+}}\rho(y_2,t)|\omega_0|(\dd y)\to 0\quad \textrm{as}\ t\to 0.
\end{align*}\\
It remains to show that 

\begin{align*}
\lim_{t\to 0}III(t)=-\int_{\mathbb{R}}u^{1}_{0}(x_1,0)\varphi(x_1,0)\dd x_1.   \tag{4.6}
\end{align*}\\
We set 

\begin{align*}
\eta(x_1,t)=2\int_{0}^{\infty}\Gamma_0(x_2,t)\varphi(x_1,x_2)\dd x_2,\quad x_1\in \mathbb{R}.  \tag{4.7}
\end{align*}\\
Since $\varphi\in C_0(\overline{\mathbb{R}^{2}_{+}})$, we have $\lim_{|x_1|\to\infty}\varphi(x_1,x_2)=0$ for each $x_2\geq 0$. Hence, $\lim_{|x_1|\to\infty} \eta(x_1,t)=0$ and $\eta(\cdot,t)\in C_{0}(\mathbb{R})$. It follows that 

\begin{align*}
\eta(x_1,t)
&=2\int_{0}^{\infty}\frac{1}{(4\pi t)^{1/2}}e^{-|x_2|^{2}/4t}\varphi(x_1,x_2)\dd x_2 \\
&=\frac{2}{\pi^{1/2}}\int_{0}^{\infty}e^{-|z_2|^{2}}\varphi(x_1,2t^{1/2}z_2)\dd z_2
\to \varphi(x_1,0)\quad \textrm{uniformly for}\ x_1\in \mathbb{R}\ \textrm{as}\ t\to 0.
\end{align*}\\
Thus, $\eta(\cdot,t)\to \varphi(\cdot,0)$ on $C_0(\mathbb{R})$ as $t\to0$. By (4.7), 

\begin{align*}
III(t)
=-2\int_{\mathbb{R}^{2}_{+}}\Gamma_0(x_2,t)\left(e^{t\partial_1^{2}}u^{1}_{0}(\cdot,0) \right)(x_1,t)\varphi(x_1,x_2)\dd x
=-\int_{\mathbb{R}}\left(e^{t\partial_1^{2}}u^{1}_{0}(\cdot,0) \right)(x_1,t)\eta(x_1,t)\dd x_1.
\end{align*}\\
Since $e^{t\partial_1^{2}}u^{1}_{0}(\cdot,0)\to u^{1}_{0}(\cdot,0)$ on $L^{1}(\mathbb{R})$, sending $t\to 0$ yields (4.6). We proved (4.3). Since $K\tilde{\omega}_0=K\omega_0$ by Lemma 2.3 (ii), $T(t)\tilde{\omega}_0=-\nabla^{\perp}\cdot S(t)K\tilde{\omega}_0=-\nabla^{\perp}\cdot S(t)K\omega_0=T(t)\omega_0$. Thus (4.4) holds.
\end{proof}

\vspace{15pt}

%remark3.2
\begin{rem}
The kernel estimate (3.21) also implies

\begin{align*}
||\partial_t^{s}\partial_x^{k}T(t)\omega_0||_{q}
\leq \frac{C}{t^{1/r-1/q+|k|/2+s}} ||\omega_0||_{r},\quad t>0,\ \omega_0\in L^{r},\ 1< r<2,\  r\leq q\leq \infty,  \tag{4.8}
\end{align*}\\
satisfying $u^{1}_0(x_1,0)=0$. The trace is understood in $L^{s,r}(\mathbb{R})$, $1/(2s)=1/r-1/2$, by $(2.7)_2$.
\end{rem}

\vspace{15pt}

\subsection{A convergence of the $L^{q}$-norm}

The function $t^{1-1/q}T(t)\omega_0$ is bounded in $L^{q}$ by (4.1). Furthermore, we have:

\vspace{15pt}

%lem3.3
\begin{lem}
For continuous measures $\omega_0\in M$, 

\begin{align*}
\lim_{t\to 0}t^{1-1/q}||T(t)\omega_0||_{q}=0,\quad 1<q\leq \infty. \tag{4.9}
\end{align*}
\end{lem}

\vspace{5pt}

\begin{proof}
We set $u_0=K\omega_0$ and $\omega=T(t)\omega_0$. Since the trace of $v=S(t)u_0$ vanishes for $t>0$, taking the trace for $v=K\omega$ implies 

\begin{align*}
\int_{\mathbb{R}^{2}_{+}}\frac{y_2}{|x_1-y_1|^{2}+y_2^{2} }\omega(y,t)\dd y=0,\quad x_1\in \mathbb{R},\ t>0.
\end{align*}\\
Thus applying (4.8) yields

\begin{align*}
t^{1-1/q}||T(t)\omega_0||_{q}
= t^{1-1/q}||T(t/2)T(t/2)\omega_0||_{q}
\lesssim t^{1-1/r}||T(t/2)\omega_0||_{r},\quad 1<r<2,\ r\leq q\leq \infty.
\end{align*}\\
It suffices to show (4.9) for $q\in (1,2)$. By (1.12), 

\begin{align*}
T(t)\omega_0=e^{t\Delta_N}\omega_0-(e^{t\Delta_N}-e^{t\Delta_D})(H\partial_2-\partial_1)\partial_1(-\Delta_D)^{-1}\omega_0-2\Gamma_0(x_2,t)e^{t\partial^{2}_{1}}u^{1}_{0}(\cdot,0).
\end{align*}\\
We shall show that 

\begin{align*}
\lim_{t\to 0}t^{1-1/r}||e^{t\Delta_N}\omega_0||_{q}=0.  \tag{4.10}
\end{align*}\\
We set a measure $\overline{\omega}_0$ on $\mathbb{R}^{2}$ by 

\begin{align*}
\overline{\omega}_{0}(E)=\omega_0(E\cap \overline{\mathbb{R}^{2}_{+}})
\end{align*}\\
for Borel sets $E\subset \mathbb{R}^{2}$. By the reflection ${\overline{\omega}}_0^{*}(E)=\overline{\omega}_{0}(E^{*})$ and $E^{*}=\{x^{*}|\ x\in E \}$, we define $\omega_0^{\textrm{even}}=\overline{\omega}_0+\overline{\omega}_0^{*}$. By changing the variable, we see that 

\begin{align*}
e^{t\Delta_{N}}\omega_0
=\int_{\mathbb{R}^{2}_{+}}\left(\Gamma(x-y,t)+\Gamma(x-y^{*},t) \right)\omega_0(\dd y)
=\int_{\mathbb{R}^{2}}\Gamma(x-y,t)\omega_0^{\textrm{even}}(\dd y)=e^{t\Delta}\omega_0^{\textrm{even}}.
\end{align*}\\
Since $\omega_0$ is continuous, so is $\omega_0^{\textrm{even}}$. Since $\lim_{t \to 0}t^{1-1/q}||e^{t\Delta}\omega_0^{\textrm{even}}||_q=0$ for the continuous measure $\omega_0^{\textrm{even}}$ \cite[Lemma 4.4]{GMO}, (4.10) follows. Since $u^{1}_{0}(\cdot,0)\in L^{1}(\mathbb{R})$ by $(2.7)_1$, we have 

\begin{align*}
t^{1-1/q}||\Gamma_0(\cdot,t)e^{t\partial_1^{2}}u^{1}_{0}(\cdot,0)||_{L^{q}(\mathbb{R}^{2}_{+})}
\lesssim t^{1/2(1-1/q)}||e^{t\partial_1^{2}}u^{1}_{0}(\cdot,0)||_{L^{q}(\mathbb{R})}\to 0\quad \textrm{as}\ t\to0.
\end{align*}\\
It remains to show that 

\begin{align*}
\lim_{t\to 0}t^{1-1/q}||(e^{t\Delta_N}-e^{t\Delta_D})(H\partial_2-\partial_1)\partial_1(-\Delta_D)^{-1}\omega_0||_{q}=0.  \tag{4.11}
\end{align*}\\
We may assume that $\omega_0\geq 0$ by the Jordan decomposition. Since 

\begin{align*}
\omega_0(\overline{\mathbb{R}^{2}_{+}})=\omega_0\left(\cup_{n\geq 1}\left\{y_2\geq   1/n\right\}\right)=\lim_{n\to\infty}\omega_0(\{y_2\geq 1/n\}),
\end{align*}\\
by continuity of the measure from below, for $\varepsilon>0$ there exists $\delta>0$ such that 

\begin{align*}
\omega_0(\{y_2< \delta\})=\omega_0(\overline{\mathbb{R}^{2}_{+}})-\omega_0(\{y_2\geq \delta\})\leq \varepsilon.
\end{align*}\\
We use the kernel representation (3.23). Since $\tilde{W}(x,y,t)$ satisfies (3.21), we estimate 

\begin{align*}
\int_{\mathbb{R}^{2}_{+}}|\tilde{W}(x,y,t)|^{q}\dd x\lesssim 
\int_{\mathbb{R}^{2}_{+}}\left|\frac{e^{-cx_2^{2}/t}}{|x-y^{*}|^{2}+t}  \right|^{q}\dd x 
\lesssim \int_{0}^{\infty}\frac{e^{-cqx_2^{2}/t}}{\left((x_2+y_2)^{2}+t\right)^{q-1/2}}\dd x_2
\lesssim \frac{t^{1/2}}{(y_2^{2}+t)^{q-1/2}}.
\end{align*}\\
It follows from (3.23) that

\begin{align*}
t^{1-1/q}||(e^{t\Delta_N}-e^{t\Delta_D})(H\partial_2-\partial_1)\partial_1(-\Delta_D)^{-1}\omega_0||_{q}
&\lesssim t^{1-1/q}\int_{\mathbb{R}^{2}_{+}}\left(\int_{\mathbb{R}^{2}_{+}}|\tilde{W}(x,y,t)|^{q}\dd x \right)^{1/q}|\omega_0|(\dd y) \\
&\lesssim \int_{\mathbb{R}^{2}_{+}}\frac{  t^{1-1/2q}  }{(y_2^{2}+t)^{1-1/2q}} |\omega_0|(\dd y) \\
&=\int_{\{y_2\geq \delta\}}+\int_{\{y_2< \delta\}} \\
&\lesssim  \frac{t^{1-1/2q}}{\delta^{2-1/q}}||\omega_0||_{M}+\varepsilon.
\end{align*}\\
Hence, $\overline{\lim}_{t\to 0}t^{1-1/q}||(e^{t\Delta_N}-e^{t\Delta_D})(H\partial_2-\partial_1)\partial_1(-\Delta_D)^{-1}\omega_0||_{q}\leq C\varepsilon$ for some constant $C>0$. Since $\varepsilon>0$ is arbitrary, (4.11) holds. The proof is complete.
\end{proof}

\vspace{15pt}

\subsection{Continuity in $L^{1}$}

We prove the continuity of $T(t)\omega_0$ in $L^{1}$ for $\omega_0\in L^{1}$ satisfying $u^{1}_{0 }(x_1,0)=0$.

\vspace{15pt}

%prop4.4
\begin{prop}
For $\omega_0\in L^{1}$ satisfying $u^{1}_{0}(x_1,0)=0$ for $u_0=K\omega_0$, there exists a sequence $\{u_{0,m}\}$ such that 

\begin{equation*}
\begin{aligned}
&\omega_{0,m}=-\nabla^{\perp}\cdot u_{0,m}\in L^{1},\\
&\textrm{spt}\ \omega_{0,m}\subset \mathbb{R}^{2}_{+},\\
&\omega_{0,m}\to \omega_0\quad \textrm{on}\ L^{1}\quad \textrm{as}\ m\to\infty.
\end{aligned}
\tag{4.12}
\end{equation*}
\end{prop}

\vspace{5pt}

\begin{proof}
Let $\overline{u}_{0}$ and $\overline{\omega}_{0}$ be the zero extensions of $u_0$ and $\omega_0$ to $\mathbb{R}^{2}$. Since $u_{0}^{1}(x_1,0)=0$, $\overline{\omega}_0=-\nabla^{\perp}\cdot \overline{u}_0$ in the sense of distribution. We set $u_{0,m}(x_1,x_2)=\overline{u}_{0}(x_1,x_2-1/m)$. Then, $\omega_{0,m}=-\nabla^{\perp} \cdot u_{0,m}$ satisfies the desired property.
\end{proof}

\vspace{15pt}

%lem 4.5
\begin{lem}
For $\omega_0\in L^{1}$ satisfying $u^{1}_{0}(x_1,0)=0$ for $u_0=K\omega_0$, $T(t)\omega_0\to \omega_0$ on $L^{1}$ as $t\to 0$.
\end{lem}
 
\vspace{5pt}

\begin{proof}
By (1.12), 

\begin{align*}
T(t)\omega_0=e^{t\Delta_N}\omega_0-(e^{t\Delta_N}-e^{t\Delta_D})(H\partial_2-\partial_1)\partial_1(-\Delta_D)^{-1}\omega_0. 
\end{align*}\\
Since $e^{t\Delta_N}\omega_0\to \omega_0$ on $L^{1}$, it suffices to show 

\begin{align*}
\lim_{t\to 0}||(e^{t\Delta_N}-e^{t\Delta_D})(H\partial_2-\partial_1)\partial_1(-\Delta_D)^{-1}\omega_0||_{1}=0.  \tag{4.13}
\end{align*}\\
We first show (4.13) under the additional assumption $\textrm{spt}\ {\omega_0}\subset \mathbb{R}^{2}_{+}$. We take $\delta>0$ such that $\textrm{spt}\ {\omega_0}\subset \{x_2>\delta\}$. It follows from  (3.23) and (3.21) that

\begin{align*}
||(e^{t\Delta_N}-e^{t\Delta_D})(H\partial_2-\partial_1)\partial_1(-\Delta_D)^{-1}\omega_0||_{1}
&\lesssim
\int_{\mathbb{R}^{2}_{+}}\left(\int_{\mathbb{R}^{2}_{+}}\frac{e^{-c x_2^{2}/t }}{(|x-y^{*}|^{2}+t)}|\omega_0(y)|\dd y\right)\dd x\\
&\lesssim \int_{0}^{\infty}\left(\int_{\mathbb{R}^{2}_{+}}\frac{e^{-c x_2^{2}/t }}{((x_2+y_2)^{2}+t)^{1/2}}|\omega_0(y)|\dd y\right)\dd x_2\\
&\lesssim \frac{t^{1/2}}{\delta}||\omega_0||_{1}\to 0\quad \textrm{as}\ t\to 0.
\end{align*}\\
Thus, (4.13) holds. If $\textrm{spt}\ \omega_0\subset \overline{\mathbb{R}^{2}_{+}}
$, we take a sequence $\{u_{0,m}\}$ satisfying (4.12) and estimate 

\begin{align*}
||T(t)\omega_0-\omega_0||_{1}
&\leq ||T(t)(\omega_0-\omega_{0,m})||_{1}
+||T(t)\omega_{0,m}-\omega_{0,m}||_{1}
+||\omega_{0,m}-\omega_0||_{1} \\
&\leq C||\omega_{0,m}-\omega_0||_{1}+||T(t)\omega_{0,m}-\omega_{0,m}||_{1}.
\end{align*}\\
Since $\textrm{spt}\ \omega_{0,m}\subset \mathbb{R}^{2}_{+}$,

\begin{align*}
\varlimsup_{m\to\infty}||T(t)\omega_0-\omega_0||_{1}\leq C||\omega_0-\omega_{0,m}||_{1}.
\end{align*}\\
Since the right-hand side tends to zero as $m\to\infty$ by (4.12), the desired result follows. 
\end{proof}

\vspace{15pt}

\subsection{The asymptotic formula} 

We prove the asymptotic formula (1.22). 

%prop4.6
\begin{prop}
\begin{align*}
\lim_{t\to\infty}||e^{t\Delta }\mu ||_{1}=0,\quad \mu\in L^{1}(\mathbb{R}^{2}),\ \int_{\mathbb{R}^{2}}\mu(x)\dd x=0.  \tag{4.14}
\end{align*}
\end{prop}

\vspace{5pt}

\begin{proof}
The assertion is well known (e.g., \cite[Lemma 3.3 (i)]{Miyakawa96}). Since the total mass of $\mu$ is zero,

\begin{align*}
(e^{t\Delta }\mu)(x)=\int_{\mathbb{R}^{2}}\left(\Gamma(x-y,t)-\Gamma(x,t)\right)\mu(y)\dd y
=t^{-1}\int_{\mathbb{R}^{2}}\left(\Gamma\left((x-y)/t^{1/2},1\right)-\Gamma\left(x/t^{1/2},1\right)\right)\mu(y)\dd y.
\end{align*}\\
Integrating $e^{t\Delta}\mu$ by $x=t^{1/2}z$ and applying the dominated convergence theorem yield

\begin{align*}
||e^{t\Delta}\mu||_{1}
&\leq \int_{\mathbb{R}^{2}}t^{-1}\left|\int_{\mathbb{R}^{2}}
\Gamma\left((x-y)/t^{1/2},1\right)-\Gamma\left(x/t^{1/2},1\right)\mu(y) \dd y\right|\dd x \\
&=\int_{\mathbb{R}^{2}}\left|\int_{\mathbb{R}^{2}}
\Gamma\left(z-y/t^{1/2},1\right)-\Gamma\left(z,1\right)\mu(y) \dd y\right|\dd z \\
&\leq \int_{\mathbb{R}^{2}}\left(\int_{\mathbb{R}^{2}}\left|
\Gamma\left(z-y/t^{1/2},1\right)-\Gamma\left(z,1\right)\right|\dd z \right)|\mu(y)|\dd y \to 0\quad \textrm{as}\ t\to\infty.
\end{align*}
\end{proof}

\vspace{15pt}

%thm
\begin{thm}
\begin{align*}
\lim_{t\to\infty}t^{1-1/q}||T(t)\omega_0||_{q}=0,\quad \omega_0\in M,\ 1\leq q\leq \infty.  \tag{4.15}
\end{align*}
\end{thm}

\vspace{5pt}

\begin{proof}
By (4.1) and (4.2), we estimate 

\begin{align*}
t^{1-1/q}||T(t)\omega_0||_{q}=t^{1-1/q}||T(t/2)T(t/2)\omega_0||_{q}
\lesssim ||T(t/2)\omega_0||_{1}.
\end{align*}\\
It suffices to show (4.15) for $q=1$. Since $T(t)\omega_0\in L^{1}$ and $S(t)u_0=0$ on $\{x_2=0 \}$ for $t>0$, we may assume that $\omega_0\in L^{1}$ and $u^{1}_{0}(x_1,0)=0$ for $u_0=K\omega_0$, i.e., 

\begin{align*}
0=u^{1}_{0}(x_1,0)=\frac{1}{\pi}\int_{\mathbb{R}^{2}_{+}}\frac{y_2}{|x_1-y_1|^{2}+y_2^{2}}\omega_0(y)\dd y,\quad x_1\in \mathbb{R}.
\end{align*}\\
By integrating $u^{1}_{0}(x_1,0)$ by the $x_1$-variable, we have

\begin{align*}
\int_{\mathbb{R}^{2}_{+}}\omega_0(y)\dd y=0.   \tag{4.16}
\end{align*}\\
By (1.12),

\begin{align*}
T(t)\omega_0=e^{t\Delta_N}\omega_0-(e^{t\Delta_N}-e^{t\Delta_D})(H\partial_2-\partial_1)\partial_1(-\Delta_D)^{-1}\omega_0.
\end{align*}\\
Let $\omega_{0}^{\textrm{even}}$ be the even extension of $\omega_0$ defined by (4.5). By (4.16), $\omega_{0}^{\textrm{even}}\in L^{1}(\mathbb{R}^{2})$ and $\int_{\mathbb{R}^{2}}\omega_{0}^{\textrm{even}}\dd x=0$. Since $e^{t\Delta_{N}}\omega_0=e^{t\Delta}\omega_{0}^{\textrm{even}}$, it follows from (4.14) that 

\begin{align*}
\lim_{t\to \infty}||e^{t\Delta_{N}}\omega_0||_{1}=0.  \tag{4.17}
\end{align*}\\
We set $\omega=T(t)\omega_0$ and consider the scaling 

\begin{align*}
\omega_{\lambda}(x,t)=\lambda^{2}\omega(\lambda x.\lambda^{2}t),\quad \lambda>0.
\end{align*}\\
Since $||\omega_{\lambda}||_{1}(t)=||\omega||_{1}(\lambda^{2}t)$, observe that 

\begin{align*}
\lim_{t\to \infty}||T(t)\omega_0||_{1}=\lim_{\lambda \to \infty}||\omega||_{1}(\lambda^{2})
=\lim_{\lambda \to \infty}||\omega_{\lambda}||_{1}(1).
\end{align*}\\
By the scaling property of the kernel (3.19), $\omega_{\lambda}=T(t)\omega_{0,\lambda}$. Hence, 

\begin{align*}
\omega_{\lambda}(\cdot,1)=e^{\Delta_N}\omega_{0,\lambda}-(e^{\Delta_N}-e^{\Delta_D})(H\partial_2-\partial_1)\partial_1(-\Delta_D)^{-1}\omega_{0,\lambda}.
\end{align*}\\
By (4.17), 

\begin{align*}
\lim_{\lambda\to\infty}||e^{\Delta_N}\omega_{0,\lambda}||_{1}
=\lim_{\lambda\to\infty}||e^{\lambda^{2}\Delta_N}\omega_{0}||_{1}
=0.
\end{align*}\\
It suffices to show that 

\begin{align*}
\lim_{\lambda\to \infty}||(e^{\Delta_N}-e^{\Delta_D})(H\partial_2-\partial_1)\partial_1(-\Delta_D)^{-1}\omega_{0,\lambda}||_{1}=0.  \tag{4.18}
\end{align*}\\
We use the kernel representation (3.23). We shall show that

\begin{align*}
\int_{\mathbb{R}^{2}_{+}}|\tilde{W}(x,y,1)|\dd x\leq C\eta(y_2),\quad  y\in \mathbb{R}^{2}_{+},  \tag{4.19}
\end{align*}\\
for some constant $C>0$ and

\begin{align*}
\eta(y_2)=\int_{0}^{y_2}\rho \Gamma_0(\rho,1)\dd \rho+y_2\int_{y_2}^{\infty} \Gamma_0(\rho,1)\dd \rho,\quad y_2\geq 0.
\end{align*}\\
The function $\eta$ is bounded continuous in $[0,\infty)$ and satisfies $\eta(0)=0$. The convergence (4.18) follows from (4.19) since 

\begin{align*}
||(e^{\Delta_N}-e^{\Delta_D})(H\partial_2-\partial_1)\partial_1(-\Delta_D)^{-1}\omega_{0,\lambda}||_{1}
&=\int_{\mathbb{R}^{2}_{+}}\left|\int_{\mathbb{R}^{2}_{+}}\tilde{W}(x,y,1)\omega_{0,\lambda }(y)\dd y\right|\dd x  \\
&\lesssim  \int_{\mathbb{R}^{2}_{+}}\eta(y_2)\lambda^{2} |\omega_0(\lambda y)|\dd y \\
&=\int_{\mathbb{R}^{2}_{+}}\eta(z_2/\lambda) |\omega_0(z)|\dd z
\to 0\quad \textrm{as}\ \lambda\to\infty.
\end{align*}\\
To prove (4.19), we use the shorthand notation $\tilde{W}(x,y)=\tilde{W}(x,y,1)$ and $\Gamma(x)=\Gamma(x,1)$. Since $\partial_{x_1}E(x)=-2^{-1}Q_{x_2}(x_1)$ by (3.6), 

\begin{align*}
\tilde{W}(x,y)
=4\int_{0}^{y_2}\int_{\mathbb{R}}\Gamma(x-z^{*})\partial_{z_1}^{2}E(z-y)\dd z
&=4\partial_{x_1}\int_{0}^{y_2}\int_{\mathbb{R}}\Gamma_0(x_1-z_1)\Gamma_0(x_2+z_2)\partial_{z_1}E(z-y)\dd z \\
&=-2\partial_{x_1}\int_{0}^{y_2}\int_{\mathbb{R}}\Gamma_0(x_1-z_1)\Gamma_0(x_2+z_2)Q_{y_2-z_2}(z_1-y_1)\dd z \\
&=-2\partial_{x_1}\int_{0}^{y_2}\Gamma_0(x_2+z_2)(Q_{y_2-z_2}*_{1}\Gamma_0)(x_1-y_1)\dd z_2.
\end{align*}\\
Since $Q_{s}*_1\Gamma_0=e^{sA}H\Gamma_0$ and $A=-H\partial_1$ by (3.4) and (3.3), 

\begin{align*}
\tilde{W}(x,y)
=-2\partial_{x_1}\int_{0}^{y_2}\Gamma_0(x_2+z_2)e^{(y_2-z_2)A}H\Gamma_0\dd z_2
=2\int_{0}^{y_2}\Gamma_0(x_2+z_2)e^{(y_2-z_2)A}A\Gamma_0\dd z_2.
\end{align*}\\
Since $A\Gamma_0\in L^{1}(\mathbb{R})$ by (3.8), we estimate

\begin{align*}
\int_{\mathbb{R}^{2}_{+}}|\tilde{W}(x,y)|\dd x
&\lesssim \int_{0}^{\infty}\left(\int_{0}^{y_2}\Gamma_0(x_2+z_2)||e^{(y_2-z_2)A}A\Gamma_0||_{L^{1}(\mathbb{R})}\dd z_2\right) \dd x_2\\
&\lesssim \int_{0}^{y_2}\left(\int_{0}^{\infty} \Gamma_0(x_2+z_2)  \dd x_2\right)\dd z_2
=\int_{0}^{y_2}\left(\int_{z_2}^{\infty} \Gamma_0(\rho)  \dd \rho\right)\dd z_2.  
\end{align*}\\
By changing the order of the integrals,

\begin{align*}
\int_{0}^{y_2}\left(\int_{z_2}^{\infty}\Gamma_0(\rho)\dd \rho\right) \dd z_2 
&=\int_{0}^{y_2}\left(\int_{0}^{\rho}\Gamma_0(\rho)\dd z_2 \right)\dd \rho
+\int_{y_2}^{\infty}\left(\int_{0}^{y_2}\Gamma_0(\rho)\dd z_2 \right)\dd \rho \\
&=\int_{0}^{y_2}\rho \Gamma_0(\rho)\dd \rho
+y_2\int_{y_2}^{\infty}\Gamma_0(\rho)\dd \rho
=\eta(y_2).
\end{align*}\\
Thus (4.19) holds. The proof is now complete. 
\end{proof}

\vspace{15pt}

\section{Applications to the Navier-Stokes flow}

\vspace{15pt}

We prove Theorems 1.1 and 1.2. 

\vspace{15pt}

\subsection{Local well-posedness}

For $u_0\in L^{2,\infty}_{\sigma}$ satisfying $\omega_0\in M$, we set

\begin{align*}
\omega_0=\omega_{0,\textrm{pp}}+\omega_{0,\textrm{cont}}.  \tag{5.1}
\end{align*}\\
Since $t^{1-1/q}T(t)\omega_{0,\textrm{cont}}\to 0$ in $L^{q}$ for $q\in (1,2)$ as $t\to 0$ by (4.9), $(2.4)_2$ and (2.12) imply that $u_{0,\textrm{cont}}=K\omega_{0,\textrm{cont}}$ satisfies 

\begin{align*}
\lim_{t\to 0}t^{1/2-1/p}(||S(t)u_{0,\textrm{cont}}||_{p}
+t^{1/2}||\nabla S(t)u_{0,\textrm{cont}}||_{p})=0,\quad 2<p<\infty.
\end{align*}\\
Thus by (2.13) and $(2.4)_1$, there exists a constant $C_{1}>0$ such that 

\begin{align*}
\varlimsup_{t\to 0}t^{1/2-1/p}(||S(t)u_0||_{p}+t^{1/2}||\nabla S(t)u_0||_{p} )
\leq C_{1}||\omega_{0,\textrm{pp}}||_{M},\quad 2< p< \infty.  \tag{5.2}
\end{align*}\\
We set a sequence $\{u_j\}$ by 

\begin{equation*}
\begin{aligned}
u_{j+1}&=u_1-\int_{0}^{t}S(t-s)\mathbb{P}(u_j\cdot \nabla u_{j})\dd s,  \\
u_1&=S(t)u_0.
\end{aligned}
\tag{5.3}
\end{equation*}\\
By taking the rotation to (5.3),

\begin{equation*}
\begin{aligned}
\omega_{j+1}&=\omega_1+\int_{0}^{t}\nabla^{\perp}\cdot S(t-s)\mathbb{P}(\omega_ju_j^{\perp})\dd s,\quad u_{j}=K\omega_j,   \\ 
\omega_1&=T(t)\omega_0.
\end{aligned}
\tag{5.4}
\end{equation*}\\
For $T>0$, set 

\begin{equation*}
\begin{aligned}
N_j&=\sup_{0< t\leq T}t^{1-1/q}||\omega_j||_{q},\quad 1\leq q<2,\\
L_j&=\sup_{0<  t\leq T}t^{1/2-1/p}\left(||u_j||_{p}+t^{1/2}||\nabla u_j||_{p}   \right),\quad 2<p<\infty.\\
\end{aligned}
\tag{5.5}
\end{equation*}\\
By the Sobolev inequality $||u_j||_{\infty}\lesssim ||u_j||_{p}^{1-2/p}||\nabla u_j||_{p}^{2/p}$, 

\begin{align*}
\sup_{0< t\leq T}t^{1/2}||u_j||_{\infty}\leq C L_j.  \tag{5.6}
\end{align*}\\

%prop5.1
\begin{prop}
\begin{align*}
L_{j+1}&\leq L_1+C_2L_{j}^{2},   \tag{5.7}\\
N_{j+1}&\leq N_1+C_2L_{j}N_{j} \tag{5.8},
\end{align*}\\
with some constant $C_2>0$.
\end{prop}

\vspace{5pt}

\begin{proof}
We set $r=p/2$. Applying (2.12) and (2.15) implies 

\begin{align*}
||u_{j+1}||_{p}
&\leq ||u_1||_{p}+\int_{0}^{t}||S((t-s)/2)S((t-s)/2)\mathbb{P}\D\ u_ju_j||_{p}\dd s\\
&\leq ||u_1||_{p}+C\int_{0}^{t}\frac{\dd s}{(t-s)^{1/r-1/p+1/2}}||u_j||_{2r}^{2}\dd s\\
&\leq ||u_1||_{p}+CL_{j}^{2}\int_{0}^{t}\frac{\dd s}{(t-s)^{1/p+1/2}s^{1-2/p} }
= ||u_1||_{p}+\frac{C'}{t^{1/2-1/p}}L_j^{2}.
\end{align*}\\
We estimate 

\begin{align*}
||\nabla u_{j+1}||_{p}
\leq ||\nabla u_1||_{p}
+\int_{0}^{t/2}||\nabla S(t-s)\mathbb{P}\D\ (u_ju_j)||_{p}\dd s
+\int_{t/2}^{t}||\nabla S(t-s)\mathbb{P}u_j\cdot \nabla u_j||_{p}\dd s.
\end{align*}\\
Applying (2.12) and (2.15) yields

\begin{align*}
\int_{0}^{t/2}||\nabla S((t-s)/2)S((t-s)/2)\mathbb{P}\D\ (u_ju_j)||_{p}\dd s
\leq CL_j^{2}\int_{0}^{t/2}\frac{\dd s}{(t-s)^{1+1/p} s^{1-2/p}} 
= \frac{C'}{t^{1-1/p}}L_j^{2},
\end{align*}
\begin{align*}
\int_{t/2}^{t}||\nabla S(t-s)\mathbb{P}u_j\cdot \nabla u_j||_{p}\dd s
&\leq C\int_{t/2}^{t}\frac{1}{(t-s)^{1/2+1/r-1/p}}||u_j\cdot \nabla u_j||_{r} \dd s \\
&\leq CL_j^{2}\int_{t/2}^{t}\frac{1}{(t-s)^{1/2+1/p}s^{3/2-2/p}} \dd s =\frac{C'}{t^{1-1/p}}L_j^{2}.
\end{align*}\\We obtained

\begin{align*}
||\nabla u_{j+1}||_{p}
\leq ||\nabla u_1||_{p}+\frac{C}{t^{1-1/p}}L_j^{2}.
\end{align*}\\
Thus (5.7) holds. By (2.14), (2.16), (5.5) and (5.6), we estimate

\begin{align*}
||\omega_{j+1}||_{q}
&\leq ||\omega_1||_{q}+\int_{0}^{t}||\nabla^{\perp}\cdot S(t-s)\mathbb{P}(\omega_ju^{\perp}_{j})||_{q}\dd s \\
&\leq ||\omega_1||_{q}+C\int_{0}^{t}\frac{1}{(t-s)^{1/2}}||(\omega_ju^{\perp}_{j})||_{q}\dd s \\
&\leq ||\omega_1||_{q}+C'L_jN_j
\int_{0}^{t}\frac{\dd s}{(t-s)^{1/2}s^{3/2-1/q}}  
= ||\omega_1||_{q}+\frac{C''}{t^{1-1/q}}L_jN_j.
\end{align*}\\
Thus (5.8) holds.
\end{proof}

\vspace{10pt}

%prop5.2
\begin{prop}
There exits a constant $\delta_1>0$ such that for $\omega_0\in M$ satisfying $||\omega_{0,\textrm{pp}}||_{M}\leq \delta_1$, there exists $T>0$ and a unique $(\omega,u)$ satisfying (1.8), (1.9) and

\begin{align*}
&\omega\in BC_{w}([0,T]; M), \tag{5.9} \\
&u\in BC_{w}([0,T]; L^{2,\infty}),  \tag{5.10} \\
&t^{1-1/q}\omega \in BC((0,T]; L^{q}),\quad 1<q<2,  \tag{5.11}\\
&t^{1/2-1/p}u,\ t^{1-1/p}\nabla u  \in BC((0,T]; L^{p}),\quad 2<p<\infty. \tag{5.12} 
\end{align*}
\end{prop}

\vspace{5pt}

\begin{proof}
It follows from (5.2) that

\begin{align*}
\varlimsup_{T\to 0}L_1\leq C_{1}||\omega_{0,\textrm{pp}}||_{M}. \tag{5.13}
\end{align*}\\
We take $\delta_1= (8C_{1}C_2)^{-1}$. By (5.13), $\overline{\lim}_{T\to 0}L_1\leq (8C_2)^{-1}$. We take $T>0$ so that $L_1\leq (4C_2)^{-1}$. Then, (5.7) yields

\begin{align*}
L_{j+1}\leq 2L_1,\quad j\geq 0.    \tag{5.14}
\end{align*}\\
Since (5.8) and (5.14) imply

\begin{align*}
N_{j+1}
\leq N_1+C_2L_jN_{j}
\leq N_1+2C_2L_1N_j
\leq N_1+\frac{1}{2}N_j,
\end{align*}\\
we have

\begin{align*}
N_{j+1}\leq 2N_1,\quad   j\geq 0.  \tag{5.15}
\end{align*}\\
Thus $(\omega_j,u_j)$ satisfies 

\begin{align*}
t^{1-1/q}\omega_j\in BC((0,T]; L^{q}),\quad 1<q<2,  \tag{5.16}\\
t^{1/2-1/p}u_j\in BC((0,T]; L^{p}),\quad 2<p<\infty.  \tag{5.17}
\end{align*}\\
We show that 

\begin{align*}
&\omega_{j}\in BC_{w}([0,T]; M),    \tag{5.18}\\
&u_j\in BC_{w}([0,T]; L^{2,\infty}).   \tag{5.19}
\end{align*}\\
Since $\omega_j$ is bounded on $M$ by (5.14), $u_j=K\omega_j$ is bounded on $L^{2,\infty}$ by $(2.4)_1$. We shall show the weak-star continuity at $t=0$. The function $T(t)\omega_0$ is vaguely continuous on $M$ at $t=0$ by Lemma 4.1. We take an arbitrary $\varphi\in C_c^{\infty}(\overline{\mathbb{R}^{2}_{+}
})$. Let $(\cdot,\cdot)$ denote the paring for $M$ and $C_0$. It follows from (5.5), (5.6) and (2.12) that 

\begin{align*}
\left|\int_{0}^{t}(\nabla^{\perp}\cdot S(t-s)\mathbb{P}(\omega_ju^{\perp}_{j}),\varphi )\dd s \right|
&=\left|\int_{0}^{t}\left(\omega_{j} u^{\perp}_{j},S(t-s)\mathbb{P}\nabla^{\perp}\cdot \varphi   \right)\dd s\right|  \\
&\leq \int_{0}^{t}||\omega_j||_{1}(s)||u_j||_{\infty}(s)||S(t-s)\mathbb{P}\nabla^{\perp}\cdot \varphi ||_{\infty}\dd s\\
&\leq N_jL_j \int_{0}^{t}\frac{1}{s^{1/2}}||S(t-s)\mathbb{P}\nabla^{\perp}\cdot \varphi ||_{\infty}\dd s\\
&\lesssim N_jL_j\sup_{0<\rho\leq  t}\rho^{1/2}||S(\rho)\mathbb{P}\nabla^{\perp}\cdot \varphi ||_{\infty}\to 0\quad \textrm{as}\ t\to 0.
\end{align*}\\
Since $C_{c}^{\infty}(\overline{\mathbb{R}^{2}_{+}})$ is dense in $C_0(\overline{\mathbb{R}^{2}_{+}})$, $\omega_j$ is vaguely continuous on $M$ at $t=0$. This proves (5.18).

We prove (5.19). The function $S(t)u_0$ is weakly-star continuous on $L^{2,\infty}$. We take an arbitrary $\varphi\in C^{\infty}_{c}(\mathbb{R}^{2}_{+})$. Applying (5.5), (5.6) and (2.12) implies 

\begin{align*}
\left|\int_{0}^{t} (S(t-s)\mathbb{P}(\omega_j u^{\perp}_j),\varphi )\dd s \right|
&=\left|\int_{0}^{t} (\omega_j u^{\perp}_j,S(t-s)\mathbb{P} \varphi   )\dd s\right| \\
&\leq N_jL_j\int_{0}^{t}\frac{1}{s^{1/2}}||S(t-s)\mathbb{P}\varphi||_{\infty}\dd s\\
&\lesssim N_jL_j\sup_{0<\rho \leq t}\rho^{1/2}||S(\rho)\mathbb{P}\varphi||_{\infty}\to 0\quad \textrm{as}\ t\to 0.
\end{align*}\\
Since $C^{\infty}_{c}(\mathbb{R}^{2}_{+})$ is dense in $L^{2,1}(\mathbb{R}^{2}_{+})$, $u_j$ is weakly-star continuous on $L^{2,\infty}$. We proved (5.19).

We estimate $\omega_{j+1}-\omega_{j}$, $u_{j+1}-u_j$ and obtain 

\begin{equation*}
\begin{aligned}
\lim_{j\to \infty}\sup_{0< t\leq T}\left(t^{1-1/q}||\omega_{j+1}-\omega_{j}||_{q}+t^{1/2-1/p}||u_{j+1}-u_{j} ||_{p}+t^{1-1/p}||\nabla (u_{j+1}-u_{j})||_{p}  \right)=0,\\
1<q<2,\ 2<p<\infty.   
\end{aligned}
\tag{5.20}
\end{equation*}\\
Since $t^{1-1/q}\omega_j$ converges in $BC((0,T]; L^{q})$ for $q\in (1,2)$ and $t^{1/2-1/p} u_j , t^{1-1/p}\nabla u_j$ converge in $BC((0,T]; L^{p})$ for $p\in (2,\infty)$, respectively, the limit $(\omega,u)$ satisfies (5.11) and (5.12). Sending $j\to\infty$ to (5.3) and (5.4) implies (1.8) and (1.9). The weak-star continuity (5.9) and (5.10) follows in the same way as (5.18) and (5.19).

The uniqueness follows by estimating the difference of two solutions $w=u-\tilde{u}$. By (5.14), (5.13) and $\delta_1=(8C_1C_2)^{-1}$, the constructed solution satisfies

\begin{align*}
\sup_{0< s\leq T}s^{1/2-1/p}||u||_{p}\leq 2C_1\delta_1.
\end{align*}\\
Since $w$ satisfies 

\begin{align*}
w=-\int_{0}^{t}S(t-s)\mathbb{P}\D\ (wu+\tilde{u}w)\dd s,
\end{align*}\\
in the same way as the proof of Proposition 5.1, we estimate

\begin{align*}
||w||_{p}
&\leq \int_{0}^{t}||S(t-s)\mathbb{P}\D\ (wu+\tilde{u}w)||_{p}\dd s\\
&\leq \frac{C_2}{t^{1/2-1/p}}\left(\sup_{0< s\leq t}s^{1/2-1/p}||w||_{p} \right)\left\{ \left(\sup_{0< s\leq t}s^{1/2-1/p}||u||_{p}\right)+\left(\sup_{0< s\leq t}s^{1/2-1/p}||\tilde{u}||_{p}\right)  \right\} \\
&\leq \frac{C_2}{t^{1/2-1/p}}\left(\sup_{0< s\leq t}s^{1/2-1/p}||w||_{p} \right)4C_1\delta_1\\
&\leq \frac{1}{2t^{1/2-1/p}}\left(\sup_{0< s\leq t}s^{1/2-1/p}||w||_{p} \right).
\end{align*}\\
Thus, $w\equiv 0$. 
\end{proof}

\vspace{15pt}

%prop5.4
\begin{prop}
\begin{align*}
&t^{1-1/q}\omega \in BC((0,T]; L^{q}),\quad 2\leq q\leq \infty,  \tag{5.21}\\
&t^{1/2-1/p}u,\ t\nabla u\in BC((0,T]; L^{\infty}). \tag{5.22} 
\end{align*}
\end{prop}

\vspace{5pt}

\begin{proof}
The property (5.21) follows from (5.22) and (5.11). By (5.12) and the Sobolev embedding, $t^{1/2}u\in BC((0,T]; L^{\infty})$. We estimate 

\begin{align*}
||\nabla (u-S(t)u_0)||_{\infty}
&\leq \int_{0}^{t/2}||\nabla S(t-s)\mathbb{P}\D\ (uu)||_{\infty}\dd s
+\int_{t/2}^{t}||\nabla S(t-s)\mathbb{P} (u\cdot \nabla u)||_{\infty}\dd s.
\end{align*}\\
Since $t^{1/2-1/2r}u, t^{1-1/2r}\nabla u\in BC((0,T]; L^{2r})$ for $r\in (1,\infty)$ by (5.12), it follows from (2.12) and (2.15) that

\begin{align*}
\int_{0}^{t/2}||\nabla S(t-s)\mathbb{P}\D\ (uu)||_{\infty}\dd s
&=\int_{0}^{t/2}||\nabla S((t-s)/2)S((t-s)/2)\mathbb{P}\D\ (uu)||_{\infty}\dd s \\
&\lesssim \int_{0}^{t/2}\frac{1}{(t-s)^{1/2+1/r}}||S((t-s)/2)\mathbb{P}\D\ (uu)||_{r}\dd s\\
&\lesssim \int_{0}^{t/2}\frac{1}{(t-s)^{1+1/r}}||u||_{2r}^{2}\dd s\\
&\lesssim \left(\sup_{0< s\leq t}s^{1-1/2r}||u||_{2r}\right) \int_{0}^{t/2}\frac{\dd s}{(t-s)^{1+1/r}s^{1-1/r}}\\
&\lesssim \left(\sup_{0< s\leq t}s^{1-1/2r}||u||_{2r}\right)\frac{1}{t}.
\end{align*}\\
Applying (2.12) yields 

\begin{align*}
\int_{t/2}^{t}||\nabla S(t-s)\mathbb{P} (u\cdot \nabla u)||_{\infty}\dd s 
&=\int_{t/2}^{t}||\nabla S((t-s)/2)S((t-s)/2)\mathbb{P}u\cdot \nabla u||_{\infty}\dd s\\
&\lesssim \int_{t/2}^{t}\frac{1}{(t-s)^{1/2+1/2r}}||u||_{\infty}||\nabla u||_{2r}\dd s\\
&\lesssim \left(\sup_{0< s\leq t}s^{1/2}||u||_{\infty}\right)\left(\sup_{0< s\leq t}s^{1-1/2r}||\nabla u||_{2r}\right)\int_{t/2}^{t}\frac{\dd s}{(t-s)^{1/2+1/2r}s^{3/2-1/2r}} \\
&\lesssim \left(\sup_{0< s\leq t}s^{1/2}||u||_{\infty}\right)\left(\sup_{0< s\leq t}s^{1-1/2r}||\nabla u||_{2r}\right)
\frac{1}{t}.
\end{align*}\\
Since $t\nabla S(t)u_0\in BC((0,\infty); L^{\infty})$ by (2.13), $t\nabla u\in BC((0,T]; L^{\infty})$. We proved (5.22).
\end{proof}

\vspace{15pt}

\begin{proof}[Proof of Theorem 1.1 (i)]
If $\omega_0$ is continuous (i.e., $\omega_{0,\textrm{pp}}\equiv 0$ ), (4.9) and (5.2) yield   

\begin{align*}
&\lim_{t\to 0}t^{1-1/q}||T(t)\omega_0||_{q}=0,\quad 1<q<2,\\
&\lim_{t\to 0}t^{1/2-1/p}( ||S(t)u_0||_{p}+t^{1/2}||\nabla S(t)u_0||_{p})=0,\quad 2<p<\infty.
\end{align*}\\
Thus $L_1\to 0$ and $N_1\to 0$ as $T\to0$. By (5.14) and (5.15), the sequence $(\omega_j,u_j)$ in the proof of Proposition 5.2 satisfy

\begin{align*}
&\lim_{t\to 0}t^{1-1/q}||\omega_j||_{q}=0,\quad 1<q<2,\\
&\lim_{t\to 0}t^{1/2-1/p}( ||u_j||_{p}+t^{1/2}||\nabla u_j||_{p})=0,\quad 2<p<\infty.
\end{align*}\\
Since $t^{1-1/q}\omega_j$ and $t^{1/2-1/p}u_j, t^{1-1/p}\nabla u_j$ converge in $BC([0,T]; L^{q})$ and $BC([0,T]; L^{p})$, the limit $(\omega,u)$ satisfies 

\begin{equation*}
\begin{aligned}
&\lim_{t\to 0}t^{1-1/q}||\omega||_{q}=0,\quad 1<q<2,\\
&\lim_{t\to 0}t^{1/2-1/p}( ||u||_{p}+t^{1/2}||\nabla u||_{p})=0,\quad 2<p<\infty.
\end{aligned}
\tag{5.23}
\end{equation*}\\
From the proof of Proposition 5.3, (5.23) holds also for $q\in [2,\infty]$ and $p\in [2,\infty]$.

If $\omega_0\in L^{1}$ and $u^{1}_{0}(x_1,0)=0$, $T(t)\omega_0\in BC([0,T]; L^{1})$ by Lemma 4.5. Applying (2.16) yields

\begin{align*}
||\omega-T(t)\omega_0||_{1}
&\leq \int_{0}^{t}||\nabla^{\perp}\cdot S(t-s)\mathbb{P}(\omega u^{\perp})||_{1}\dd s \\
&\lesssim \int_{0}^{t}\frac{1}{(t-s)^{1/2}}||\omega||_{1} ||u||_{\infty}\dd s\\
&\lesssim \left(\sup_{0<s\leq t}||\omega||_1\right)\left(\sup_{0<s\leq t}s^{1/2}||u||_{\infty}\right)\to 0\quad \textrm{as}\ t\to 0.
\end{align*}\\
Thus, $\omega\in BC([0,T]; L^{1})$ follows. By $(2.4)_1$, $u\in BC([0,T]; L^{2,\infty})$. The proof is complete.
\end{proof}

\vspace{15pt}

\subsection{Global well-posedness}
It remains to show Theorems 1.1 (ii) and 1.2.

\vspace{15pt}

\begin{proof}[Proof of Theorem 1.1 (ii)]
For $u_0\in L^{2,\infty}_{\sigma}$ satisfying $\omega_0\in M$, (2.13), $(2.4)_1$ and (4.1) yield

\begin{align*}
\sup_{0<t<\infty}t^{1/2-1/p}\left(||S(t)u_0||_{p}+t^{1/2}||\nabla S(t)u_0||_{p}\right)
&\leq C_3||\omega_0||_{M},\quad 2<p<\infty, \\
\sup_{0<t<\infty}t^{1-1/q}||T(t)\omega_0||_{q}
&\leq C_3||\omega_0||_{M},\quad  1\leq q<2. 
\end{align*}\\
We set a sequence $(\omega_j,u_j)$ by (5.3), (5.4) and take $T=\infty$ in (5.5). Then, (5.6)-(5.8) holds. We assume $||\omega_0||_{M}\leq \delta_2$ for $\delta_2= (8C_2C_3)^{-1}$. Then, $L_1\leq (4C_2)^{-1}$ and (5.14)-(5.20) holds for $T=\infty$. Thus the limit $(\omega,u)$ satisfies (1.8), (1.9), (5.9)-(5.12), (5.21) and (5.22). The uniqueness follows in the same way as the proof of Proposition 5.2.
\end{proof}

\vspace{15pt}

\begin{proof}[Proof of Theorem 1.2]
Let $(\omega,u)$ be a local-in-time solution in $[0,t_0]$ constructed in Theorem 1.1 (i). We take an arbitrary $T>0$. Since $u(\cdot,t_0)\in L^{\infty}$, $\nabla u(\cdot,t_0)\in L^{2}$ and $u(\cdot,t_0)=0$ on $\{x_2=0\}$, by \cite[Remark 6.4]{A7}, $u$ is extendable to a global-in-time solution in $[t_0,T]$ satisfying 

\begin{align*}
u,\ (t-t_0)^{1/2}\nabla u\in BC_w([t_0,T]; L^{\infty}),\quad \nabla u\in BC_{w}([t_0,T]; L^{2}),
\end{align*}\\
and 

\begin{align*}
u(t)=S(t-t_0)u(t_0)-\int_{t_0}^{t}S(t-s)\mathbb{P}(u\cdot \nabla u)\dd s,\quad t_0\leq t\leq T. \tag{5.25}
\end{align*}\\
Taking the rotation yields 

\begin{align*}
\omega(t)=T(t-t_0)\omega(t_0)+\int_{t_0}^{t}\nabla^{\perp}\cdot S(t-s)\mathbb{P}(\omega u^{\perp})\dd s,\quad t_0\leq t\leq T.  \tag{5.26}
\end{align*}\\
Since (1.18) and (1.20) hold in $(0,t_0]$, we have 

\begin{align*}
t^{1/2}u,\ t\nabla u\in BC((0,T]; L^{\infty}),\quad  
t^{1/2}\nabla u\in BC((0,T]; L^{2}).     \tag{5.27}
\end{align*}\\
By the H\"oler's inequality,

\begin{align*}
t^{1-1/p}\nabla u\in BC((0,T]; L^{p}),\quad 2<p\leq \infty.     \tag{5.28}
\end{align*}\\
By (5.25) and (5.28), (1.20) follows. We take $q\in (1,2)$. It follows from (2.14), (5.28) and (1.20) that 

\begin{align*}
||\omega||_q(t)\lesssim ||\omega||_{q}(t_0)+\int_{t_0}^{t}\frac{1}{(t-s)^{1/2}}||\omega||_{2q}(s)||u||_{2q}(s)\dd s 
\lesssim ||\omega||_{q}(t_0)+\int_{t_0}^{t}\frac{1}{(t-s)^{1/2}s^{3/2-1/q}}\dd s.
\end{align*}\\
Thus (1.19) holds. By (1.20) and (1.19), 

\begin{align*}
||\omega u^{\perp}||_{1}\leq ||\omega||_{p'} ||u||_{p}\lesssim \frac{1}{t^{1/2}},\quad 0<t\leq T.
\end{align*}\\
By (5.26) and (2.16), (1.17) follows. Since (1.17) implies (1.18) by $(2.4)_1$, we proved (1.17)-(1.20) for $T>0$.
\end{proof}

\vspace{15pt}

\appendix

\section{A solution formula under the trace zero condition}

\vspace{15pt}

We show that the formula (1.4) gives a solution to (3.24). 

\vspace{15pt}

%lem A.1
\begin{thm}
Let $q\in (1,2)$. Let $\omega_0\in L^{q}$ satisfy $u^{1}_{0}(x_1,0)=0$ for $u_0=K\omega_0$. Then, $\omega=T_0(t)\omega_0$ satisfies (3.24) and $T_0(t)\omega_0\to \omega_0$ in $L^{q}$ as $t\to 0$.
\end{thm}

\vspace{5pt}

\begin{proof}
By (1.4), 

\begin{equation*}
\begin{aligned}
W_0(x,y,t)
&=\Gamma(x-y,t)-\Gamma(x-y^{*},t)+2(H\partial_{1}-\partial_{2})\partial_2E*\Gamma (x-y^{*},t). \\
T(t)\omega_0&=e^{t\Delta_D}\omega_0+2(H\partial_1-\partial_2)\int_{\mathbb{R}^{2}_{+}}(\partial_2E*\Gamma)(x-y^{*},t)\omega_0(y)\dd y.
\end{aligned}
\tag{A.1}
\end{equation*}\\
Since $(\partial_t -\Delta_{x})W_0(x,y,t)=0$, $\omega=T(t)\omega_0$ satisfies the heat equation. Since $A=-H\partial_1$, it follows that 

\begin{align*}
(\partial_2-A)(H\partial_1-\partial_2)\partial_2E*\Gamma(z,t)
=-(\partial_2-A)(\partial_2+A)\partial_2 E*\Gamma(z,t)
=\partial_2\Gamma(z,t).
\end{align*}\\
Multiplying $\partial_2-A$ by $W_0(x,y,t)$ yields 

\begin{align*}
(\partial_2-A)W_0(x,y,t)
&=(\partial_2-A)(\Gamma(x-y,t)-\Gamma(x-y^{*},t))
+2(\partial_2-A)(H\partial_1-\partial_2)\partial_2E*\Gamma(x-y^{*},t) \\
&=(\partial_2-A)(\Gamma(x-y,t)-\Gamma(x-y^{*},t))
+2\partial_2\Gamma(x-y^{*},t)\\
&=\partial_2(\Gamma(x-y,t)+\Gamma(x-y^{*},t))
-A(\Gamma(x-y,t)-\Gamma(x-y^{*},t))\\
&\to 0\quad \textrm{as}\ x_2\to0.
\end{align*}\\
Thus, $(\partial_2-A)\omega=0$ on $\{x_2=0\}$.  We shall show the convergence to initial data. We observe that 

\begin{align*}
\int_{\mathbb{R}^{2}_{+}}\partial_{x_2}E(x-y^{*})\omega_0(y)\dd y=0,\quad x\in \mathbb{R}^{2}_{+}.     \tag{A.2}
\end{align*}\\
The left-hand side belongs to $L^{p}$, $1/p=1/q-1/2$, by $(2.4)_2$ and is harmonic in $\mathbb{R}^{2}_{+}$, vanishing on $\{x_2=0\}$ by $u^{1}_{0}(x_1,0)=0$. Thus, (A.2) follows from the Liouville theorem. With the operators $(-\Delta_D)^{-1}$ and 

\begin{align*}
(-\Delta_N)^{-1}\omega_0=\int_{\mathbb{R}^{2}_{+}}(E(x-y)+E(x-y^{*}))\omega_0(y)\dd y,
\end{align*}\\
(A.2) is represented as

\begin{align*}
\partial_2\left(  (-\Delta_N)^{-1}-(-\Delta_D)^{-1} \right)\omega_0=0.  \tag{A.3}
\end{align*}\\
Since 

\begin{align*}
(\partial_iE*\Gamma)(z,t)=\int_{t}^{\infty}\partial_i\Gamma(z,s)\dd s, \quad i=1,2, \tag{A.4}
\end{align*}\\
we have

\begin{align*}
2\int_{\mathbb{R}^{2}_{+}}(\partial_2E*\Gamma)(x-y^{*},t)\omega_0(y)\dd y
&=2\int_{\mathbb{R}^{2}_{+}}\left(\int_{t}^{\infty}\partial_2 \Gamma(x-y^{*},s)\dd s\right)\omega_0(y)\dd y\\
&=\int_{t}^{\infty}\partial_{2}(e^{s\Delta_N}-e^{s\Delta_D})\omega_0\dd s \\
&=\partial_{2}\left((-\Delta_{N})^{-1}e^{t\Delta_N}-(-\Delta_{D})^{-1}e^{t\Delta_D}\right)\omega_0\\
&=\partial_{2}\left((-\Delta_{N})^{-1}(e^{t\Delta_N}-I)-(-\Delta_{D})^{-1}(e^{t\Delta_D}-I)  \right)\omega_0.
\end{align*}\\
Since $H$ and $\nabla^{2}(-\Delta_N)^{-1}, \nabla^{2}(-\Delta_D)^{-1}$ are bounded on $L^{q}(\mathbb{R})$ and $L^{q}(\mathbb{R}^{2}_{+})$, it follows from (A.1) that 

\begin{align*}
||T(t)\omega_0-\omega_0||_{q}
\lesssim ||e^{t\Delta_N}\omega_0-\omega_0||_{q}
+||e^{t\Delta_D}\omega_0-\omega_0||_{q}\to 0\quad \textrm{as}\ t\to 0.
\end{align*}\\
The proof is complete.
\end{proof}

\vspace{15pt}

%remsA.2
\begin{rems}
(i) The formula 

\begin{equation*}
\begin{aligned}
W_0(x,y,t)
&=\Gamma(x-y,t)+\Gamma(x-y^{*},t)-2A(\partial_2+A)E*\Gamma(x-y^{*},t),\\
T_0(t)\omega_0
&=e^{t\Delta_N}\omega_0-2A(\partial_2+A)\int_{\mathbb{R}^{2}_{+}}E*\Gamma(x-y^{*},t)\omega_0(y)\dd y,
\end{aligned}
\tag{A.5}
\end{equation*}\\
is obtained in \cite[Theorem 3.1]{Maekawa13}. Since $A=-H\partial_1$ and 

\begin{align*}
(H\partial_1-\partial_2)\partial_2E=H\partial_1\partial_2E+\partial_1^{2}E-\Delta E
=-A\partial_2 E-A^{2}E +\delta_0
=-A(\partial_2+A)E+\delta_0,
\end{align*}\\
we have 

\begin{align*}
W_0(x,y,t)
&=\Gamma(x-y,t)-\Gamma(x-y^{*},t)+2(H\partial_1-\partial_2)\partial_2E*\Gamma(x-y^{*},t)\\
&=\Gamma(x-y,t)+\Gamma(x-y^{*},t)-2A(\partial_2+A)E*\Gamma(x-y^{*},t).
\end{align*}\\

\noindent
(ii) We are able to derive (A.5) from (3.24). Indeed, for a solution $\omega$ of (3.24), $(\partial_2-A)\omega$ satisfies the heat equation subject to the Dirichlet boundary condition, i.e.,   

\begin{align*}
(\partial_2-A)\omega=e^{t\Delta_D}(\partial_2-A)\omega_0.
\end{align*}\\
Since (3.13) implies 

\begin{align*}
e^{t\Delta_D}(\partial_2-A)\omega_0
=\partial_2 e^{t\Delta_N}\omega_0-Ae^{t\Delta_D}\omega_0
=(\partial_2-A)e^{t\Delta_N}\omega_0+A(e^{t\Delta_N}-e^{t\Delta_D})\omega_0,
\end{align*}\\
we have 

\begin{align*}
\Delta \omega=(\partial_2+A)(\partial_2-A)\omega
&=(\partial_2+A)e^{t\Delta_D}(\partial_2-A)\omega_0 \\
&=(\partial_2+A)(\partial_2-A)e^{t\Delta_N}\omega_0 
+(\partial_2+A)A(e^{t\Delta_N}-e^{t\Delta_D})\omega_0  \\
&=\Delta e^{t\Delta_N}\omega_0+(\partial_2+A)A(e^{t\Delta_N}-e^{t\Delta_D})\omega_0.
\end{align*}\\
Integrating $\partial_t\omega=\Delta\omega$ on $(0,t)$ yields

\begin{align*}
\omega
=e^{t\Delta_N}\omega_0
+A(\partial_2+A)\int_{0}^{t}(e^{s\Delta_N}-e^{s\Delta_D})\omega_0\dd s.
\end{align*}\\
Since (A.3) holds for $\omega_0$ satisfying $u^{1}_{0}(x_1,0)=0$ for $u_0=K\omega_0$, 

\begin{align*}
A(\partial_2+A)\int_{0}^{\infty}(e^{s\Delta_N}-e^{s\Delta_D})\omega_0\dd s
&=-H\partial_1(\partial_2-H\partial_1)((-\Delta_N)^{-1}-(-\Delta_D)^{-1} )\omega_0\\
&=-(H\partial_1\partial_2+\partial_1^{2})((-\Delta_N)^{-1}-(-\Delta_D)^{-1} )\omega_0\\
&=-(H\partial_1\partial_2-\partial_2^{2}+\Delta )((-\Delta_N)^{-1}-(-\Delta_D)^{-1} )\omega_0\\
&=-(H\partial_1-\partial_2)\partial_2((-\Delta_N)^{-1}-(-\Delta_D)^{-1} )\omega_0\\
&=0.
\end{align*}\\
Hence by (A.4), we have 

\begin{align*}
A(\partial_2+A)\int_{0}^{t}(e^{s\Delta_N}-e^{s\Delta_D})\omega_0\dd s
&=-A(\partial_2+A)\int_{t}^{\infty}(e^{s\Delta_N}-e^{s\Delta_D})\omega_0\dd s \\
&=H\partial_1(\partial_2+A)\int_{t}^{\infty}(e^{s\Delta_N}-e^{s\Delta_D})\omega_0\dd s \\
&=2H(\partial_2+A)\int_{t}^{\infty}\left(\int_{\mathbb{R}^{2}_{+}}\partial_1\Gamma(x-y^{*},t)\omega_0(y)\dd y\right)\dd s \\
&=2H(\partial_2+A)\int_{\mathbb{R}^{2}_{+}}\partial_1E*\Gamma(x-y^{*},t)\omega_0(y)\dd y \\
&=-2A(\partial_2+A)\int_{\mathbb{R}^{2}_{+}}E*\Gamma(x-y^{*},t)\omega_0(y)\dd y.
\end{align*}\\
Thus (A.5) follows.

\noindent 
(iii) The kernel $W_0(x,y,t)$ has an explicit form

\begin{align*}
W_{0}(x,y,t)=\Gamma(x-y,t)+\Gamma(x-y^{*},t)
+4\int_{\mathbb{R}^{2}_{+}}\partial_{z_1}^{2}E(z)\Gamma(x-y^{*}-z,t)\dd z.  \tag{A.6}
\end{align*}\\
Since

\begin{align*}
A(\partial_2+A)E(z)=-\partial_1(H\partial_2+\partial_1)E(z)
=
\begin{cases}
\quad -2\partial_1^{2}E(z)&\quad z_2>0,\\
\quad 0  & \quad z_2<0,
\end{cases}
\end{align*}\\
by (3.7) and 

\begin{align*}
-2A(\partial_2+A)E*\Gamma(z,t)=-2\int_{\mathbb{R}^{2}}A(\partial_2+A)E(w)\Gamma(z-w)\dd w
=4\int_{\mathbb{R}^{2}_{+}}\partial_1^{2}E(w)\Gamma(z-w,t)\dd w,
\end{align*}\\
(A.6) follows from (A.5). Since the kernel (A.6) is not integrable for the $x$-variable in $\mathbb{R}^{2}_{+}$, the formula (A.5) does not give a $L^{1}$-bound for $T_0(t)\omega_0$ and $\omega_0\in M$. However, if $u^{1}_{0}(x_1,0)=0$ for $u_0=K\omega_0$, $T_0(t)\omega_0$ agrees with $T(t)\omega_0$ by the following Theorem A.3. Hence $T_0(t)\omega_0$ belongs to $L^{1}$ for such $\omega_0$ by (4.1).
\end{rems}

\vspace{15pt}

%thmA.2
\begin{thm}
For $\omega_0\in M$ (resp. $\omega_0\in  L^{q}$, $q\in (1,2)$) satisfying $u^{1}_{0}(x_1,0)=0$ for $u_0=K\omega_0$, 

\begin{align*}
T_0(t)\omega_0=T(t)\omega_0.   \tag{A.7}
\end{align*}
\end{thm}

\vspace{5pt}

\begin{proof}
Since (A.2) holds by $u^{1}_{0}(x_1,0)=0$ and $H\partial_{x_2}E(x)=\partial_{x_1}E(x)$ for $x\in \mathbb{R}^{2}_{+}$ by (3.7), multiplying $H$ by (A.2) implies 

\begin{align*}
\int_{\mathbb{R}^{2}_{+}}\partial_{x_1}E(x-y^{*})\omega_0(\dd y)=0,\quad x\in \mathbb{R}^{2}_{+}.   \tag{A.8}
\end{align*}\\
We use the indicator function 

\begin{align*}
\chi(z_2)
=\begin{cases}
\ 1,& z_2\geq 0,  \\
\ 0,& z_2<0,
\end{cases}
\end{align*}\\
to observe that

\begin{align*}
\int_{\mathbb{R}^{2}_{+}}\partial_{z_1}^{2}E(z)\Gamma(x-y^{*}-z,t)\dd z
&=\partial_{1}^{2}E\chi*\Gamma (x-y^{*},t) \\
&=\int_{\mathbb{R}^{2}}\partial_{z_1}^{2}E(z-y^{*})\chi(z_2+y_2)\Gamma(x-z,t)\dd z \\
&=\int_{\mathbb{R}^{2}}\partial_{w_1}^{2}E(w^{*}-y^{*})\chi(-w_2+y_2)\Gamma(x-w^{*},t)\dd w.
\end{align*}\\
Since $\chi(-w_2+y_2)=0$ for $w_2>y_2$ and 

\begin{align*}
 \partial_{x_1}^{2}E(x)=-\frac{x_2^{2}-x_1^{2}}{2\pi |x|^{4}}
 =\partial_{x_1}^{2}E(x^{*}),
\end{align*}\\
changing the variable yeilds

\begin{align*}
&\int_{\mathbb{R}^{2}}\partial_{w_1}^{2}E(w^{*}-y^{*})\chi(-w_2+y_2)\Gamma(x-w^{*},t)\dd w \\
&=\int_{0}^{y_2}\int_{\mathbb{R}}\partial_{w_1}^{2}E(w^{*}-y^{*})\Gamma(x-w^{*},t)\dd w
+\int_{-\infty}^{0}\int_{\mathbb{R}}\partial_{w_1}^{2}E(w^{*}-y^{*})\Gamma(x-w^{*},t)\dd w \\
&=\int_{0}^{y_2}\int_{\mathbb{R}}\partial_{z_1}^{2}E(z-y)\Gamma(x-z^{*},t)\dd z
+\int_{0}^{\infty}\int_{\mathbb{R}}\partial_{z_1}^{2}E(z-y^{*})\Gamma(x-z,t)\dd z.
\end{align*}\\
We integrate the both sides by the measure $\omega_0$. It follows from (A.8) that

\begin{align*}
\int_{\mathbb{R}^{2}_{+}}\left(\int_{\mathbb{R}^{2}_{+}}\partial_{z_1}^{2}E(z)\Gamma(x-y^{*}-z,t)\dd z\right)\omega_0(\dd y) 
=\int_{\mathbb{R}^{2}}\left(\int_{0}^{y_2}\int_{\mathbb{R}}\partial_{z_1}^{2}E(z-y)\Gamma(x-z^{*},t)\dd z\right)\omega_0(\dd y).
\end{align*}\\
It follows from (A.6) that 

\begin{align*}
T_0(t)\omega_0
&=\int_{\mathbb{R}^{2}_{+}}W_0(x,y,t)\omega_0(\dd y) \\
&=e^{t\Delta_N}\omega_0
+4\int_{\mathbb{R}^{2}_{+}} \left(\int_{0}^{y_2}\int_{\mathbb{R}}\partial_{z_1}^{2}E(z-y)\Gamma(x-z^{*},t)\dd z \right)\omega_0(\dd y)\\
&=T(t)\omega_0.
\end{align*}\\
The proof is complete.
\end{proof}

\vspace{5pt}

\section*{Acknowledgements}
This work is partially supported by JSPS through the Grant-in-aid for Young Scientist (B) 17K14217, Scientific Research (B) 17H02853 and Osaka City University Strategic Research Grant 2018 for young researchers.

%ref
\bibliographystyle{abbrv}
\bibliography{ref}

\end{document}